\documentclass[12pt]{article}

\usepackage{amsmath}
\usepackage{amsthm}
\usepackage{amssymb}
\usepackage{graphics}
\usepackage{graphicx}

\newtheorem{theorem}{Theorem}
\newtheorem{lemma}[theorem]{Lemma}
\newtheorem{corrol}[theorem]{Corollary}

\def\defn#1{{\bf #1}}
\def\tabelm#1{{\begin{tabular}[t]{c} #1 \end{tabular}}}
\def\multi#1{{\begin{tabular}[t]{@{}l@{}}#1\end{tabular}}}

\def\Real{{\mathbb R}}
\def\Cmpx{{\mathbb C}}
\def\Intg{{\mathbb Z}}
\def\Natn{{\mathbb N}}
\def\cnj#1{{\overline{#1}}}

\def\sp{{\,,\quad}}
\def\sdotsc{,\!..,}

\def\Poisson(#1,#2){{\left\{#1,#2\right\}}}
\def\Set#1{{\left\{#1\right\}}}

\def\ket(#1){{\left|#1\right>}}
\def\bra(#1){{\left<#1\right|}}

\def\Aelement#1{{\boldsymbol{#1}}}
\def\Aa{\Aelement{a}}
\def\Aw{\Aelement{w}}
\def\Ax{\Aelement{x}}
\def\Ay{\Aelement{y}}
\def\Az{\Aelement{z}}
\def\Af{\Aelement{f}}
\def\Ag{\Aelement{g}}
\def\Av{\Aelement{v}}
\def\Au{\Aelement{u}}

\def\Aeps{\Aelement{\varepsilon}}
\def\Aalpha{\Aelement{\alpha}}

\def\Aq{\Aelement{q}}

\def\AF{\Aelement{F}}
\def\AC{\Aelement{C}}

\def\CC(#1){C_{#1}}

\def\Idn{{\boldsymbol{1}\hspace{-0.25em}\textup{\small I}}}

\def\Hil{{\cal H}}
\def\Lin{{\cal L}}
\def\Alg{{\cal A}}
\def\Balg{{\cal B}}
\def\Bpi{{\pi_\Balg}}
\def\man{{\cal M}}

\def\SetS{{\cal S}}

\title{A Noncommutative Geometric Analysis of a Sphere/Torus Topology
Change}

\author{J Gratus\thanks{email: j@gratus.net}}

\begin{document}
\maketitle

\begin{abstract}
A one parameter set of noncommutative complex algebras
is given. These may be considered deformation quantisation algebras.
The commutative limit of these algebras correspond to the
algebra of polynomial functions over a manifold or variety.  The
topology of the manifold or variety depends on the parameter, varying
from nothing, to a point, a sphere, a certain variety and finally a
torus. The irreducible adjoint preserving representations of the
noncommutative algebras are studied.  As well as typical
noncommutative sphere type representations and noncommutative torus
type representations, a new object is discovered and called a
Sphere-Torus.
\end{abstract}

\tableofcontents

%%%%%%%%%%%%%%%%%%%%%%%%%%%%%%%%%%%%%%%%%%%%%%%%%%%%%%%%%%%%%%%%%%%%%%

\section{Introduction}
\label{ch_intr}

In noncommutative geometry we often wish to find analogues to
topological properties of a manifold such as compactness,
connectedness and genus, which we will consider here. For matrix
geometries~\cite{Madore_book} the question of genus is tricky since
both the sphere and torus have matrix analogues.

In deformation quantisation~\cite{Sternheimer1}, 
one can simply take the commutative limit
and ask what the genus of the underlying manifold is.  In the case of
the sphere and torus, we can also take the representation of the
noncommutative algebra and compare its properties with the
representations of the noncommutative sphere and torus.

In this article we present, in section \ref{ch_alg}, a one parameter
set of deformation algebras $\Alg(R)$ for $R\in\Real$. The commutative
limit of these algebras are $C^\omega_0(\man(R))$ the commutative
algebra of complex polynomials on the manifold (or variety)
$\man(R)$. This manifold has different topologies depending on the
value of $R$. Varying from nothing, to a point, a sphere, a variety
and finally a torus. This is described in section
\ref{ch_comm}. 

In section \ref{ch_rep} we look at the finite dimensional
representations of $\Alg(R)$. These can be classified as either
$S^2$-type representations or $T^2$-type representations by comparing
them to the representations for the noncommutative torus or the
noncommutative sphere. Depending on the value of $R$ one or other of
these representations exists. What we show in this article is that
there is a region of $R$ where both types of representation
exist. This region, which we shall name the \defn{sphere-torus}, is a
purely noncommutative region, it disappears in the commutative limit.

We summarise the various representation in the conclusion, section
\ref{ch_concl}. 

%%%%%%%%%%%%%%%%%%%%%%%%%%%%%%%%%%%%%%%%%%%%%%%%%%%%%%%%%%%%%%%%%%%%%%

\subsection{Notation}
\begin{tabular}{lp{13cm}}
$\Real_+$ 
&
$\Set{t\in\Real\,|\,t>0}$
\\
$\man$
&
Manifold or variety.
\\
$C^\omega(\man)$, $C^\omega(\Real^r)$ 
&
Commutative algebra of complex
analytic functions over $\man$ or $\Real^r$.
\\
$C^\omega_0(\man)$, $C^\omega_0(\Real^r)$ 
&
Commutative algebra of complex
polynomials over $\man$ or $\Real^r$.
\\
$\Alg$, $\Balg$ 
&
Noncommutative algebras over $\Cmpx$.
\\
\multi{Unbold symbols:\\\quad $f,g,u,v,x_i,$\\\quad $F_s,w,x,y,z$} 
&
Elements of the commutative
algebras. i.e. Analytic or polynomial functions over $\man$ or
$\Real$.
\\
\multi{Bold symbols:\\\quad $\Af,\Ag,\Au,\Av,\Ax_i,$\\
\quad $\AF_s,\Aw,\Ax,\Ay,\Az$} 
&
Elements of the noncommutative algebras $\Alg$ and $\Balg$
\\
$\Aeps$ 
&
Element in the centre 
of the noncommutative algebras $\Alg$ and $\Balg$
\\
$\Hil$ 
&
Finite or infinite dimensional Hilbert space.
\\
$\Lin(\Hil,\Hil)$
&
Space of linear maps over $\Hil$. 
(Bounded or
unbounded operators.)
\\
$M_n(\Cmpx)$ 
&
Space of $n\times n$ complex matrices.
\\
$\Idn$
&
The identity in $\Lin(\Hil,\Hil)$ and the unit matrix in
$M_n(\Cmpx)$.
\\
$\ket(\theta)$
&
Bra-ket notation for vectors in $\Hil$.
\end{tabular}

%%%%%%%%%%%%%%%%%%%%%%%%%%%%%%%%%%%%%%%%%%%%%%%%%%%%%%%%%%%%%%%%%%%%%%

\subsection{A one parameter set of immersions with a sphere, 
torus and variety}
\label{ch_comm}

\newlength{\pichigh}
\setlength{\pichigh}{3.2 cm}
\setlength{\unitlength}{\pichigh}
\begin{figure}
\begin{tabular}{|c|c|c|c|}
\hline
\begin{picture}(0.8,1.1)
\put(0,0){\includegraphics[height=\pichigh]{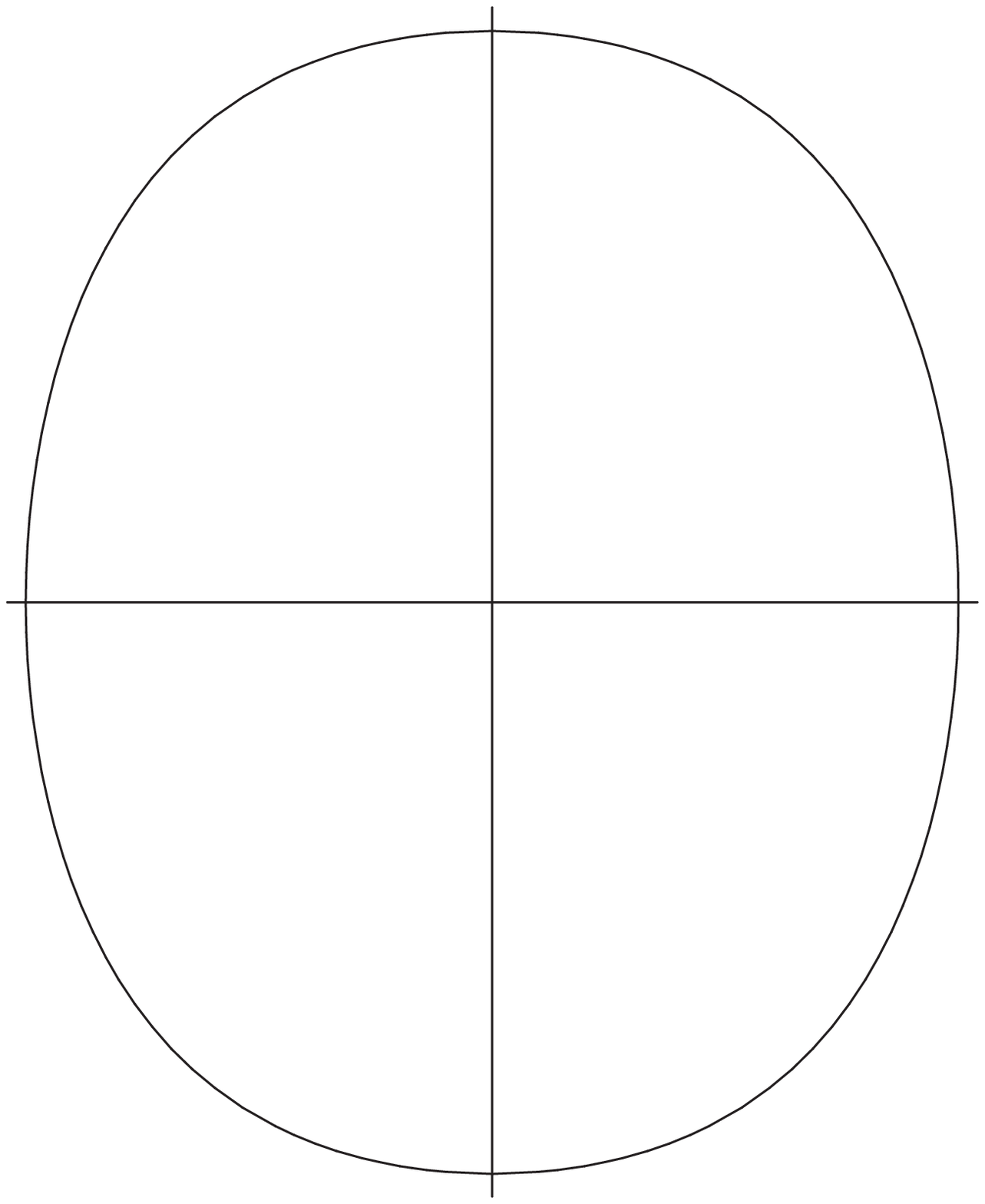}}
\put(0.4,0.8){\makebox(0,0)[rt]{$z\,$}}
\end{picture}
&
\begin{picture}(1.4,1.1)
\put(0,0){\includegraphics[height=\pichigh]{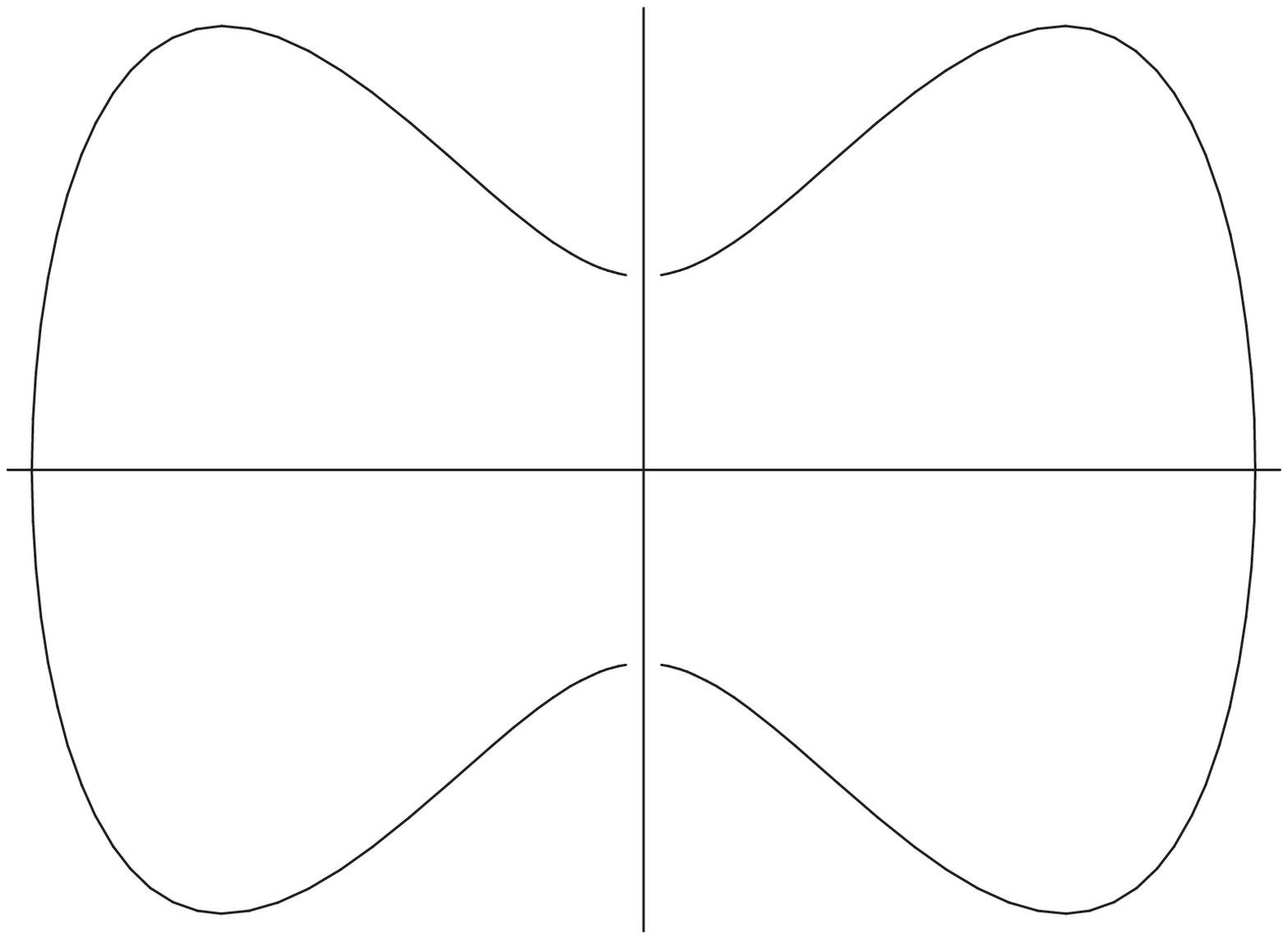}}
\put(0.7,0.9){\makebox(0,0)[rt]{$z\,$}}
\end{picture}
&
\begin{picture}(1.4,1.1)
\put(0,0){\includegraphics[height=\pichigh]{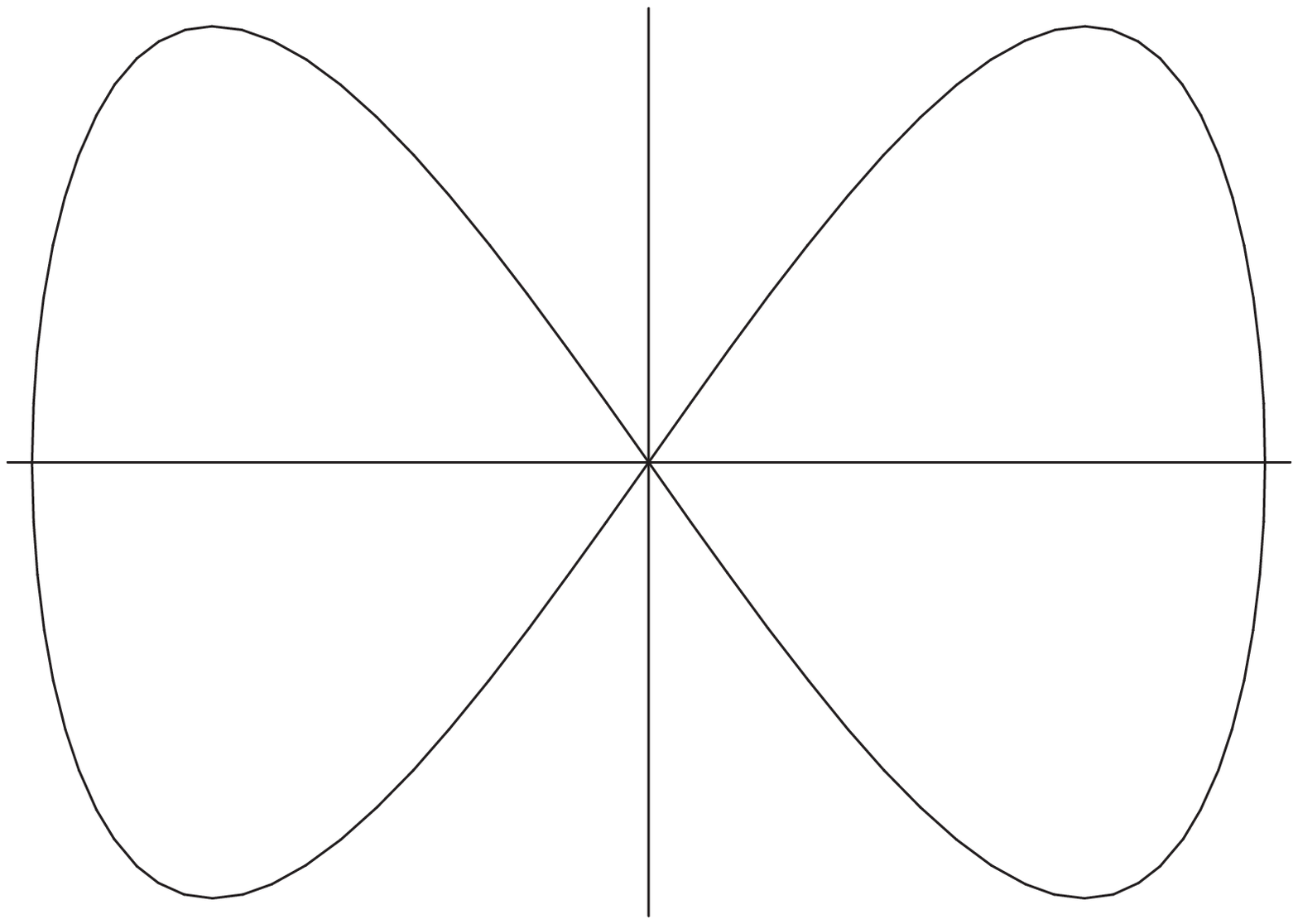}}
\put(0.7,0.9){\makebox(0,0)[rt]{$z\,$}}
\end{picture}
&
\begin{picture}(1.4,1.1)
\put(0,0){\includegraphics[height=\pichigh]{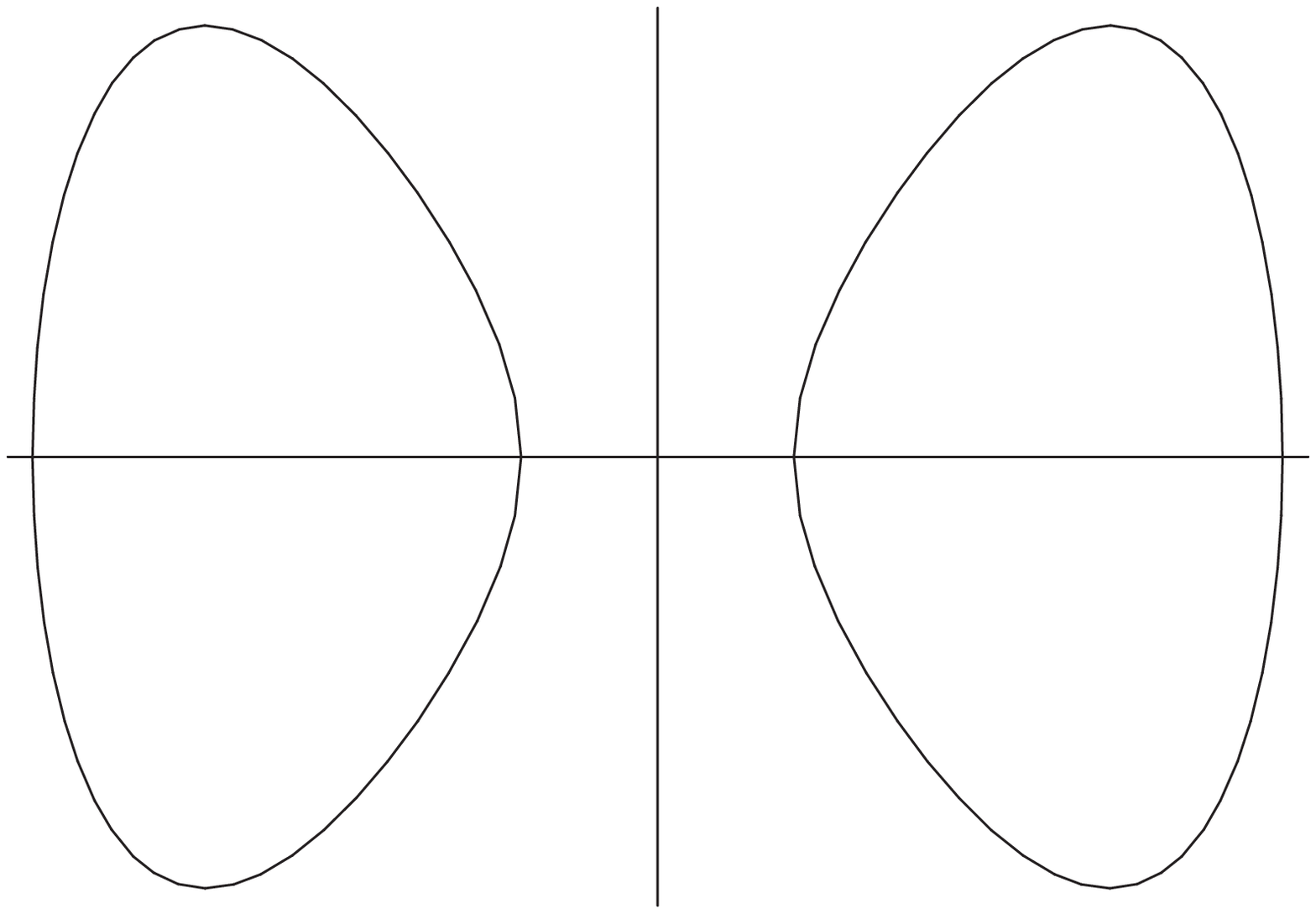}}
\put(0.7,0.9){\makebox(0,0)[rt]{$z\,$}}
\end{picture}
\\
$-1<R\le0$ 
&
$0<R<1$
&
$R=1$ 
&
$R>1$ 
\\
\hline
\end{tabular}
\caption{Shape of the slice of $\man$ setting $y=0$. The shape of
$\man$ is given by rotating slice about the $z$-axis.}
\label{fig_slices}
\end{figure}

Consider the immersions given by
%[
\begin{align}
\man(R)=\Set{(x,y,z)\in\Real^3\,\bigg|\,z^2+(x^2+y^2-R)^2=1}
\label{intr_man}
\end{align}
%]
It is obvious that for $R<-1$ there is no solution to (\ref{intr_man})
while for $R=-1$, $\man(R)$ consists of the single point at the origin
$(x,y,z)=(0,0,0)$. To picture $\man(R)$ for $R>-1$, we note that it is
axisymmetric about the $z$-axis. Therefore we can examine the shape of
$\man(R)$ by setting $y=0$ and rotating the subsequent 1 dimensional
variety about the $z$-axis. From figure \ref{fig_slices} we can see
that for $-1<R\le0$, $\man(R)$ is a convex manifold topologically
equivalent to the sphere.  For $0<R<1$, $\man(R)$ is a non convex
manifold topologically equivalent to the sphere.  For $R=1$, $\man(R)$
is not a manifold but instead a 2 dimensional variety, which is smooth
about all points except the origin $(x,y,z)=(0,0,0)$.  For $R>1$,
$\man(R)$ is a torus.

There is a Poisson structure on $\man$ is given by
%[
\begin{align}
\Poisson(x,y)=z \sp
\Poisson(z,x)=2(x^2+y^2-R)y=2wy \sp
\Poisson(y,z)=2(x^2+y^2-R)x=2wx 
\label{intr_poi}
\end{align}
%]
where $w=x^2+y^2-R$, which is consistent with (\ref{intr_man})

We can give a Darboux coordinate system $(p,q)$ such that
%[
\begin{align}
x = (R+\cos(2p))^{1/2}\cos(q)
\sp
y = -(R+\cos(2p))^{1/2}\sin(q)
\sp
z = \sin(2p)
\label{intr_darb}
\end{align}
%]
where $\Poisson(p,q)=1$. It is easy to see that (\ref{intr_man}) and
(\ref{intr_poi}) are satisfied. A necessary condition for these to be
valid is $R+\cos(2p)>0$. Thus
for the torus $R>1$ this is valid for all $p$. More
specifically we can patch coordinate systems with $0<q<2\pi$ and
$0<p<\pi$. For the variety $R=1$ we must exclude the point
$p=\tfrac12\pi$ which correspond to the point at the origin.  For
$-1<R<1$ we have
$-\tfrac12(\pi-\arccos(R))<p<\tfrac12(\pi-\arccos(R))$, and we must
exclude the two points $(x,y,z)=(0,0,\pm(1-R^2)^{1/2})$.

%%%%%%%%%%%%%%%%%%%%%%%%%%%%%%%%%%%%%%%%%%%%%%%%%%%%%%%%%%%%%%%%%%%%%%

\subsection{Brief introduction to deformation quantisation}

We limit ourselves in this article to a the deformation quantisation
of algebraic manifolds and varieties with algebraic Poisson structures.
Let $\man$ smooth $m$ dimensional Poisson manifold or variety given by
%[
\begin{equation}
\begin{aligned}
\man=\Set{(x_1\sdotsc x_r)\in\Real^r
\bigg| 
F_s(x_1\sdotsc x_r)=0
,\,
s=1\sdotsc r\!-\!m
}
\end{aligned}
\label{intr_gen_man}
\end{equation}
%]
where $F_s(x_1\sdotsc x_m)$ are polynomials. Let the Poisson structure
$\Poisson(\bullet,\bullet)$ be given by
$\Poisson(x_i,x_j)=C_{ij}(x_1\sdotsc x_m)$, where
$C_{ij}(x_1\sdotsc x_m)$ are also polynomials. Consistency implies
$\Poisson(x_i,F_s)=0$.  Let $C^\omega(\man)$ be the algebra of
complex analytic functions on $\man$ and let $C^\omega_0(\man)$ be the
subalgebra of polynomials in $(x_1\sdotsc x_m)$, this is dense in
$C^\omega(\man)$.

Let $\Balg$ be the free noncommutative
algebra generated by $\Set{\Ax_1\sdotsc \Ax_r,\Aeps}$, and define
the linear map $\Bpi:\Balg\mapsto C^\omega_0(\Real^r)$ where
$C^\omega_0(\Real^r)$ is the algebra of polynomials in
$(x_1\sdotsc x_n)$, via 
%[
\begin{align}
\Bpi(\Af+\Ag)=\Bpi(\Af)+\Bpi(\Ag)\sp
\Bpi(\Af\Ag)=\Bpi(\Af)\Bpi(\Ag)\sp
\Bpi(\lambda)=\lambda\sp
\Bpi(\Aeps)=0\sp
\Bpi(\Ax_i)=x_i
\end{align}
%]
Choose the elements $\AF_s\in\Balg$ for $s=1\sdotsc r\!-\!m$ and 
$\AC_{ij}\in\Balg$ for $i,j=1\sdotsc r$ such that 
$\Bpi(\AF_s)=F_s$ and $\Bpi(\AC_{ij})=C_{ij}$. 
We define the algebra $\Alg$ to be $\Balg$ quotiented by the noncommutative
polynomial relationships
%[
\begin{align}
[\Aeps,\Ax_i]=0\sp
\AF_s = 0 \sp
[\Ax_i,\Ax_j] = \Aeps\AC_{ij}
\label{intr_A_rels}
\end{align}
%]
for $i,j=1\sdotsc r$ and $s=1\sdotsc r-m$. 
The first equation in (\ref{intr_A_rels}) implies that $\Aeps$ is in
the centre of $\Alg$. 
We demand that $\Alg$ be an associative algebra. This imposes
restrictions on the possible choices of $\AF_s$ and 
$\AC_{ij}$, which we will not investigate here.

We can define the map
%[
\begin{align}
\pi:\Alg\mapsto C^\omega_0(\man)\,;\quad \pi(\Aeps)=0
\label{intr_pi}
\end{align}
%]
which is surjective. It is easy to see that this gives the Poisson 
structure via
%[
\begin{align}
\Poisson({\pi(\Af)},{\pi(\Ag)})=\pi\left(\frac{1}{i\Aeps}[\Af,\Ag]\right)
\end{align}
%]
Thus the following diagram commutes
%[
\begin{align}
\begin{array}{c}
\setlength{\unitlength}{1em}
\begin{picture}(17,6)
\put(0.7,5){$\Balg$}
\put(0,0){$C^\omega_0(\Real^r)$}
\put(1,4.5){\vector(0,-1){3.5}}
\put(-0.5,3){$\Bpi$}
\put(3.5,0.5){\vector(1,0){14}}
\put(2,5.5){\vector(1,0){15.5}}
\put(3,6){$/$\footnotesize 
           $\Set{[\Ax_i,\Aeps]=0,\, 
            \AF_s = 0,\, [\Ax_i,\Ax_j] = \Aeps\AC_{ij}}$}
\put(7,1){$/$\footnotesize $\Set{F_s = 0}$}
\put(18.5,5){$\Alg$}
\put(18,0){$C^\omega_0(\man)$}
\put(19,4.5){\vector(0,-1){3.5}}
\put(19.5,3){$\pi$}
\end{picture}
\end{array}
\end{align}
%]
We also demand that there is a conjugate structure $\dagger:\Alg\mapsto\Alg$
such that
%[
\begin{align}
(\Af\Ag)^\dagger=\Ag^\dagger\Af^\dagger
\sp
\lambda^\dagger=\cnj{\lambda} \text{ for } \lambda\in\Cmpx
\sp
\pi(\Af^\dagger)=\cnj{\pi(\Af)}
\sp
\Aeps^\dagger=\Aeps
\label{intr_conj}
\end{align}
%]

%%%%%%%%%%%%%%%%%%%%

We are interested in representations 
$\Psi:\Alg\mapsto\Lin(\Hil,\Hil)$ where $\Lin(\Hil,\Hil)$ is the space
of linear maps on the Hilbert space $\Hil$.
We demand that $\Psi$ is irreducible, and 
$\Psi(\Aeps)=\varepsilon\Idn$ where $\varepsilon\in\Real_+$.
If $\dim\Hil=n$ then $\Lin(\Hil,\Hil)\cong M_n(\Cmpx)$.  On the other hand if
$\dim\Hil=\infty$ then $\Lin(\Hil,\Hil)$ may contain unbounded
operators. We also demand that $\Psi$ preserve the conjugate structure
$\Psi(\Af^\dagger)=\Psi(\Af)^\dagger$ where the dagger on the right is
the Hermitian conjugate or adjoint.

We use the bra-ket notation so that if $\ket(\theta)\in\Hil$ then
$\Psi(\Af)\ket(\theta)$ the action of $\Af\in\Alg$ on $\ket(\theta)$
is written $\Af\ket(\theta)$.
Since $\Psi$ preserves the
conjugate we have
$\bra(\theta')\Af^\dagger\ket(\theta)=\cnj{\bra(\theta)\Af\ket(\theta')}$.

We sometimes want to recover the commutative structure from the matrix
algebras. This requires finding a sequence of representations
$\Psi_n:\Alg\mapsto M_n(\Cmpx)$, $\Psi_n(\Aeps)\to0$ as $n\to\infty$. 
However in general there is no canonical map $M_n(\Cmpx)\mapsto
M_{n+1}(\Cmpx)$. One exception being the noncommutative sphere
\cite{Gratus1}.

%%%%%%%%%%%%%%%%%%%%%%%%%%%%%%%%%%%%%%%%%%%%%%%%%%%

As stated above, in this article we give a shall give a one parameter
family of such algebras, whose representations can be compared to those for
the sphere and torus. Here we give a brief summary of these two
noncommutative geometries.

\vspace{1em}

\noindent
The \defn{Noncommutative Sphere} is generated by $\Ax,\Ay,\Az,\Aeps$ with 
%[
\begin{gather}
[\Ax,\Ay]=i\Aeps\Az\sp
[\Ay,\Az]=i\Aeps\Ax\sp
[\Az,\Ax]=i\Aeps\Ay\sp
\Ax^2+\Ay^2+\Az^2=1\sp
\Ax^\dagger=\Ax\sp
\Ay^\dagger=\Ay\sp
\Az^\dagger=\Az
\end{gather}
%]
The representations are the finite dimensional irreducible
representations of $SO(3)$. A basis for $\Hil$ is
$\ket(0)\sdotsc \ket(n-1)$ so that
%[
\begin{equation}
\begin{aligned}
&
\varepsilon=2(n^2-1)^{-1/2}
\sp &&
\Aa_+\ket(r)=\varepsilon(n-r-1)^{1/2}(r+1)^{1/2}\ket(r+1)
\sp\\ 
&
\Az\ket(r)=\varepsilon (r-\tfrac12(n-1))\ket(r)
\sp &&
\Aa_-\ket(r)=\varepsilon(n-r)^{1/2}r^{1/2}\ket(r-1)
\end{aligned}
\label{intr_S2_reps}
\end{equation}
%]
where $\Aa_+=\Ax+i\Ay$ and $\Aa_-=\Ax-i\Ay$ are called the ladder
elements. We note that the ladder operators, which are the
representations of the ladder elements, terminate $\Psi_n(\Aa_+)^n=0$.
We also note that $\Psi$ is unique for each $n$, and that
$\Psi_n(\Aeps)\to0$ as $n\to\infty$.

\vspace{1em}

In our language we the  \defn{Noncommutative Torus} is generated by
$\Set{\Ax_1,\Ax_2,\Ax_3,\Ax_4,\Aeps}$ with 
%[
\begin{equation}
\begin{aligned}
&\Ax_1^2+\Ax_2^2=1\sp \\
&\Ax_3^2+\Ax_4^2=1\sp \\
&[\Ax_1,\Ax_2]=0\sp \\
&[\Ax_3,\Ax_4]=0\sp 
\end{aligned}
\quad
\begin{aligned}
&[\Ax_1,\Ax_3]=-i\Aeps(\Ax_2\Ax_4+\Ax_4\Ax_2)\sp \\
&[\Ax_1,\Ax_4]=i\Aeps(\Ax_2\Ax_3+\Ax_3\Ax_2)\sp \\
&[\Ax_2,\Ax_3]=i\Aeps(\Ax_1\Ax_4+\Ax_4\Ax_1)\sp \\
&[\Ax_2,\Ax_4]=-i\Aeps(\Ax_1\Ax_3+\Ax_3\Ax_1)\sp 
\end{aligned}
\quad
\begin{aligned}
&\Ax_1^\dagger=\Ax_1\sp \\
&\Ax_2^\dagger=\Ax_2\sp \\
&\Ax_3^\dagger=\Ax_3\sp \\
&\Ax_4^\dagger=\Ax_4
\end{aligned}
\label{intr_T2_alg}
\end{equation}
%]
If we set $\Au=\Ax_1+i\Ax_2$, $\Av=\Ax_3+i\Ax_4$, and
$\Aq=(1+i\Aeps)/(1-i\Aeps)$ then we derive the usual noncommutative
torus.
%[
\begin{align}
\Au\Au^\dagger=\Au^\dagger\Au=1
\sp
\Av\Av^\dagger=\Av^\dagger\Av=1
\sp
\Au\Av=\Aq\Av\Au
\label{intr_T2_alg2}
\end{align}
%]
However as we have defined the algebra $\Alg$, the element $\Aq$ and
$\Aq^\dagger$ are not members of $\Alg$. To solve this problem we say
that $\Alg$ is generated by
$\Set{\Ax_1,\Ax_2,\Ax_3,\Ax_4,\Aeps,(1+\Aeps^2)^{-1}}$. 

There are both finite and infinite representations of the
noncommutative torus. The finite dimensional representations $\Psi_n$
have a basis $\ket(0)\sdotsc \ket(n-1)$
%[
\begin{equation}
\begin{aligned}
&\Psi_n(\Aq)=q\Idn\sp q=e^{2\pi i k/n} \sp 
\Au\ket(r)=e^{i(\beta + 2\pi r k i / n )}\ket(r) \sp
\\
&\Av\ket(r)=\ket(r+1) \sp \Av\ket(n-1)=\nu\ket(0)
\end{aligned}
\label{intr_T2_ref_fin}
\end{equation}
%]
where $n,k\in\Natn$ and where $\nu\in\Cmpx$, $|\nu|=1$ is a phase.
We impose that $n$ and $k$ are relatively prime, so that there are no
multiple eigenvalues of $\Psi(\Au)$. There also exist other more
complicated finite dimensional representations of the noncommutative
torus where $\Psi(\Au)$ has multiple eigenvalues.

We note that the ladder elements of this representation are given by
$\Av$ and $\Av^\dagger$ and that the ladder operators do not terminate
$\Psi_n(\Av)^m\ne0$ for all $n,m\in\Intg$.  To specify $\Psi_n$
completely requires giving $(n,k,\beta,\nu)$.  If $k=1$ then
$\Psi_n(\Aq)\to1$ as $n\to\infty$.

The infinite dimensional representations have a basis $\ket(r)$,
$r\in\Intg$, and are determined by the parameters
$\alpha,\beta\in\Real$ where $\alpha/2\pi$ is irrational
%[
\begin{equation}
\begin{aligned}
\Psi(\Aq)&=q\Idn\sp q=e^{i\alpha}
\sp
\Au\ket(r)=e^{i(r\alpha+\beta)}\ket(r)
\sp
\Av\ket(r)=\ket(r+1)
\end{aligned}
\label{intr_T2_ref_inf}
\end{equation}
%]
The eigenvalues of $\Psi(\Au)$ are dense on unit circle.

%%%%%%%%%%%%%%%%%%%%%%%%%%%%%%%%%%%%%%%%%%%%%%%%%%%%%%%%%%%%%%%%%%%%%%
%%%%%%%%%%%%%%%%%%%%%%%%%%%%%%%%%%%%%%%%%%%%%%%%%%%%%%%%%%%%%%%%%%%%%%

%%%%%%%%%%%%%%%%%%%%%%%%%%%%%%%%%%%%%%%%%%%%%%%%%%%%%%%%%%%%%%%%%%%%%%
%%%%%%%%%%%%%%%%%%%%%%%%%%%%%%%%%%%%%%%%%%%%%%%%%%%%%%%%%%%%%%%%%%%%%%

\section{$\Alg(R)$: The deformation algebra of polynomials on $\man(R)$}
\label{ch_alg}

In section \ref{ch_comm}, we define a one parameter set of immersions
$\man(R)\in\Real^3$. The commutative algebra of complex valued
polynomials in $(x,y,z)$ on $\man(R)$ is written $C^\omega_0(\man(R))$
and of course is dense in $C^\omega(\man(R))$. Here we give a one
parameter set of complex noncommutative algebras $\Alg(R)$ with
$R\in\Real$, which are the deformation quantisation of
$C^\omega_0(\man(R))$.  Each $\Alg(R)$ is generated by
$\Set{\Ax,\Ay,\Az,\Aeps,(1+\Aeps^2)^{-1}}$ with $\Aeps$ in the centre
of $\Alg(R)$. The reason for including $(1+\Aeps^2)^{-1}$ is similar
as for the noncommutative torus.  These are related via
%[
\begin{align}
[\Ax,\Ay] = i\Aeps\Az \sp
[\Ay,\Az] = i\Aeps(\Aw\Ax + \Ax\Aw)  \sp
[\Az,\Ax] = i\Aeps(\Aw\Ay + \Ay\Aw)
\label{alg_com}
\end{align}
%]
and
%[
\begin{align}
\Az^2+\Aw^2 = 1 
\label{alg_cass}
\end{align}
%]
where we define the element $\Aw\in\Alg(R)$ via
%[
\begin{align}
\Aw = \Ax^2 + \Ay^2 - R
\label{alg_def_w}
\end{align}
%]
The adjoint operation is given by
%[
\begin{align}
\Ax^\dagger = \Ax \sp
\Ay^\dagger = \Ay \sp
\Az^\dagger = \Az \sp
\Aw^\dagger = \Aw 
\label{alg_adj}
\end{align}
%]
It is easy to see that, assuming (\ref{alg_com}-\ref{alg_def_w})
define an associative algebra, then $\pi:\Alg(R)\mapsto
C^\omega_0(\man(R))$, where $\pi(\Ax)=x$, $\pi(\Ay)=y$,
$\pi(\Az)=z$. It also gives with the correct Poisson structure. We
shall show that these relationships define an associative algebra in
lemma \ref{thm_alg_ass} below. But first we need to define some new
elements of $\Alg(R)$ and derive some relationships which are valid if
$\Alg(R)$ is associative.

We define the ladder elements $\Aa_+,\Aa_-\in\Alg(R)$ via 
%[
\begin{align}
\Aa_+ = \Ax + i\Ay \sp 
\Aa_- = \Ax - i\Ay 
\label{alg_def_ap_am}
\end{align}
%]
In order to emphasise the circular nature of $\Aw$ and $\Az$ we shall
define the element $\Au\in\Alg(R)$ via
%[
\begin{align}
\Au = \Aw + i\Az \in\Alg(R)
\label{alg_def_u}
\end{align}
%]
and we show below that $\Au$ is unitary, i.e. $\Au^\dagger=\Au^{-1}$.

We define the ``pseudo element'' $\Aalpha$ via
%[
\begin{align}
\Aeps = \tan(\tfrac12\Aalpha)  
\end{align}
%]
and observe that $\Aalpha$ is not a member of $\Alg(R)$. However
since $\Aeps\in\Alg(R)$ and $(1+\Aeps^2)^{-1}\in\Alg(R)$ then
from the tan half angle formulae we have
$\sin(\Aalpha)=2\Aeps(1+\Aeps^2)^{-1}\in\Alg(R)$ and 
$\cos(\Aalpha)=(1-\Aeps^2)(1+\Aeps^2)^{-1}\in\Alg(R)$. Also 
$e^{\pm i n\Aalpha}=(\cos(\Aalpha)\pm i\sin(\Aalpha))^n$ so 
$e^{i n\Aalpha}\in\Alg(R)$ for all $n\in\Intg$.

%%%%%%%%%%%%%%%%%%%%%%%%%%%%%%%%%%%%%%%%%%
Before showing that the algebra $\Alg(R)$ is associative we derive
some direct consequences of these definitions.

\begin{lemma}
From the above definitions, and the assumption that $\Alg(R)$ is
associative, we have the following relationships:
%[
\begin{align}
&
[\Az,\Aw] = 0
\label{alg_com_zw}
\\
&
\Au\Au^\dagger=\Au^\dagger\Au=1
\label{alg_u_udag}
\\
&
\Aa_+\Aa_- = \Aw + R +\Aeps \Az \sp
\Aa_-\Aa_+ = \Aw + R -\Aeps \Az
\label{alg_ap_am}
\\
&
[\Aw,\Aa_+] = -\Aeps(\Az\Aa_+ + \Aa_+\Az) \sp
[\Aw,\Aa_-] = +\Aeps(\Az\Aa_- + \Aa_-\Az)
\label{alg_com_wa}
\end{align}
%]
Also
%[
\begin{equation}
\begin{aligned}
&\Az\Aa_+ = \Aa_+(\cos(\Aalpha)\Az + \sin(\Aalpha)\Aw) 
\sp                                           
\Aw\Aa_+ = \Aa_+(-\sin(\Aalpha)\Az + \cos(\Aalpha)\Aw)
\\                                            
&\Az\Aa_- = \Aa_-(\cos(\Aalpha)\Az - \sin(\Aalpha)\Aw) 
\sp                                           
\Aw\Aa_- = \Aa_-(\sin(\Aalpha)\Az + \cos(\Aalpha)\Aw)
\end{aligned}
\label{alg_normal}
\end{equation}
%]
or alternatively
%[
\begin{align}
\Au\Aa_+=\Aa_+\Au e^{i\Aalpha} \sp
\Au\Aa_-=\Aa_-\Au e^{-i\Aalpha} \sp
\Au^\dagger\Aa_+=\Aa_+\Au^\dagger e^{-i\Aalpha} \sp
\Au^\dagger\Aa_+=\Aa_-\Au^\dagger e^{i\Aalpha}
\label{alg_normal_u}
\end{align}
%]
The general element $\Af\in\Alg(R)$ can be written uniquely in the form
%[
\begin{align}
\Af = \sum_{r=0}^N\sum_{s=-N}^N \Aa_+^r\Au^s\xi_{r,s}(\Aeps) +
\sum_{r=1}^N\sum_{s=-N}^N \Aa_-^r\Au^s\xi_{-r,s}(\Aeps) 
\label{alg_gen_el}
\end{align}
%]
for some $N\in\Natn$ where for all $r,s=-N\sdotsc N$ the function
$\xi_{r,s}(t)$ is a ration function in $t$ with denominator
$(1+t^2)^m$ for some $m\in\Natn$.
\end{lemma}
%%%%%%%%%%%%%%%%%%%%%%%%%%%%%%%%%%%%%%%%

\begin{proof}
Equation (\ref{alg_ap_am}) follows automatically from
(\ref{alg_def_ap_am}),  (\ref{alg_def_w}) and the first equation 
in (\ref{alg_com}).
Equation (\ref{alg_com_wa}) follows from
%[
\begin{align*}
[\Aw,\Aa_+] 
= 
[\Aa_+\Aa_-  - R -\Aeps \Az,\Aa_+] 
=
\Aa_+[\Aa_-,\Aa_+] - \Aeps[\Az,\Aa_+]
=
-\Aeps(2\Aa_+\Az + [\Az,\Aa_+])
\end{align*}
%]
and likewise for $[\Aw,\Aa_-]$. For (\ref{alg_com_zw}) we have
%[
\begin{align*}
[\Az,\Aw] &= 
[\Az,\Aa_+\Aa_-  - R -\Aeps \Az] = 
[\Az,\Aa_+]\Aa_- + \Aa_+[\Az,\Aa_-] 
=
\Aeps ( \Aw\Aa_+\Aa_- + \Aa_+\Aw\Aa_-    
   - \Aa_+\Aw\Aa_- - \Aa_+\Aa_-\Aw) 
\\ &=
\Aeps ( \Aw\Aa_+\Aa_- - \Aa_+\Aa_-\Aw) =
-\Aeps [\Aa_+\Aa_- ,\Aw] =
-\Aeps[\Aw + R + \Aeps \Az,\Aw] =
-\Aeps^2[\Az,\Aw] 
\end{align*}
%]
hence $(1+\Aeps^2)[\Az,\Aw]=0$, so
$[\Az,\Aw]=(1+\Aeps^2)^{-1}(1+\Aeps^2)[\Az,\Aw]=0$.
The unitarity of $\Au$
(\ref{alg_u_udag}) follow from the definition of $\Au$ and the
commutativity of $\Az$ and $\Aw$.

From (\ref{alg_com}) and (\ref{alg_com_wa}) we can have 
%[
\begin{align*}
\Az\Aa_+ - \Aa_+\Az = + \Aeps \Aw\Aa_+ +\Aeps\Aa_+\Aw
\sp
-\Aeps\Az\Aa_+ - \Aeps\Aa_+\Az = \Aw\Aa_+ - \Aa_+\Aw
\end{align*}
%]
solving these as simultaneous equations and using the tan half angle
identity give the first of results of (\ref{alg_normal}). The other
identities in (\ref{alg_normal}) follow similarly. The identities in
(\ref{alg_normal_u}) are then the complex version of
(\ref{alg_normal}).

The generators all $\Alg(R)$ are all of the form (\ref{alg_gen_el}). Given
$\Af$ of this form, then $\Af\Au$, $\Af\Aeps$ and
$\Af(1-\Aeps^2)^{-1}$ are all of the form (\ref{alg_gen_el}). 
Also 
%[
\begin{align*}
\Af\Aa_+= 
\sum_{r=0}^N\sum_{s=-N}^N \Aa_+^{r+1}\Au^s e^{i s\Aalpha}\xi_{r,s}(\Aeps) +
\sum_{r=1}^N\sum_{s=-N}^N \Aa_-^{r-1}
  (\tfrac12(1+i\Aeps)\Au+\tfrac12(1-i\Aeps)\Au^{-1}+R)
   \Au^s e^{-i s\Aalpha}\xi_{-r,s}(\Aeps) 
\end{align*}
%]
So $\Af\Aa_+$ is of the form (\ref{alg_gen_el}), likewise for
$\Af\Aa_-$.

For uniqueness, from linearity, it is enough to show that if $\Af$ is
of the form (\ref{alg_gen_el}) and $\Af=0$ then $\xi_{r,s}=0$ for all
$r,s$. By multiplying $\Af$ with $(1+\Aeps^2)^M$ for sufficiently high
$M$ then we can assume $\xi_{r,s}$ are all polynomials. Let $m$ be the
largest degree of these polynomials. Now 
%[
\begin{align*}
0 = \pi(\Af) = 
\sum_{r=0}^N\sum_{s=-N}^N \pi(\Aa_+)^r\pi(\Au)^s\xi_{r,s}(0) +
\sum_{r=1}^N\sum_{s=-N}^N \pi(\Aa_-)^r\pi(\Au)^s\xi_{-r,s}(0) 
\end{align*}
%]
so by looking at the coordinate system on $\man(R)$ this implies
$\xi_{r,s}(0)=0$ for all $r,s$. Thus
$\xi_{r,s}(\Aeps)=\Aeps\hat\xi_{r,s}(\Aeps)$ where
$\hat\xi_{r,s}(\Aeps)$ are polynomials of degree $\le m-1$. Continuing
this gives $\xi_{r,s}=0$
\end{proof}

\begin{lemma}
\label{thm_alg_ass}
The relationships (\ref{alg_def_w}) to (\ref{alg_normal_u}) define
an associative algebra.
\end{lemma}
\begin{proof}
Let $\SetS$ be set of all expressions of the form
(\ref{alg_gen_el}). We define a product $\cdot$ on $\SetS$ using the
above definitions.  

We wish to show that
$\Af_1\cdot(\Af_2\cdot\Af_3)=(\Af_1\cdot\Af_2)\cdot\Af_3$ for 
$\Af_1,\Af_2,\Af_3\in\SetS$. Where the inner bracket must be written
in the form (\ref{alg_gen_el}) first.
It is sufficient to show that 
$\Af_1\cdot(\Af_2\cdot\Af_3)=(\Af_1\cdot\Af_2)\cdot\Af_3$ 
where $\Af_1,\Af_2,\Af_3$ are
from the set $\Set{\Au,\Aa_+,\Aa_-,\Aeps,(1+\Aeps^2)^{-1}}$. This is
because any expression can be constructed from these five
elements. The element $\Au^\dagger=2\Aw-\Au$ with $\Aw$ given by
(\ref{alg_def_w}).
If $\Af_1,\Af_2$ or $\Af_3$ are either $\Aeps$ or $(1+\Aeps^2)^{-1}$
then the association relation holds since $\Aeps$ commutes with all
the generators.

The remaining 27 relationships must be proved in turn,
the interesting ones are
%[
\begin{align*}
(\Au\cdot\Aa_+)\cdot\Aa_- 
&= 
\Aa_+\Au e^{i\Aalpha} \cdot\Aa_- 
= 
(\tfrac12(1-i\Aeps)\Au+\tfrac12(1+i\Aeps)\Au^{-1}+R)\Au 
\\
&=
\Au(\tfrac12(1-i\Aeps)\Au+\tfrac12(1+i\Aeps)\Au^{-1}+R)
=
\Au\cdot(\Aa_+\cdot\Aa_-) 
\end{align*}
%]
and
%[
\begin{align*}
(\Aa_+\cdot\Aa_-)\cdot\Aa_+ 
&= 
(\tfrac12(1-i\Aeps)\Au+\tfrac12(1+i\Aeps)\Au^{-1}+R)\cdot\Aa_+ 
= 
\Aa_+(\tfrac12(1-i\Aeps)\Au e^{i\Aalpha}+
      \tfrac12(1+i\Aeps)\Au^{-1} e^{i\Aalpha}+R)
\\
&=
\Aa_+(\tfrac12(1+i\Aeps)\Au +
      \tfrac12(1-i\Aeps)\Au^{-1} +R)
=
\Aa_+\cdot(\Aa_-\cdot\Aa_+) 
\end{align*}
%]

\end{proof}

%%%%%%%%%%%%%%%%%%%%%%%%%%%%%%%%%%%%%%%%%%%%%%%%%%%%%%%%%%%%%%%%%%%%%%

\section{Finite dimensional representations of $\Alg(R)$}
\label{ch_rep}

As mentioned in the introduction, we wish to find
irreducible representations over a finite Hilbert space $\Hil$ of 
$\Alg(R)$, i.e. $\Psi:\Alg(R)\mapsto \Lin(\Hil,\Hil)\cong M_n(\Cmpx)$. 
Such that $\Psi(\Aeps)=\varepsilon\Idn$ where
$\varepsilon\in\Real_+$ and such that $\Psi$
preserves the adjoint operator so that
$\Psi(\Af^\dagger)=\Psi(\Af)^\dagger$ where the dagger on the right
hand side is the Hermitian conjugate.
We note that since $\Psi(\Au)^\dagger\Psi(\Au)=\Idn$ we can
diagonalise $\Psi(\Au)$ and the eigenspaces of $\Psi(\Au)$ are
orthogonal. We shall further assume that $\Psi(\Au)$ has no multiple
eigenvalues, so the eigenspaces of $\Psi(\Au)$ are all one
dimensional.  This significantly simplifies the types of
representations.

\begin{theorem}
\label{th_reps}
Let $\Psi:\Alg(R)\mapsto\Lin(\Hil,\Hil)\cong M_n(\Cmpx)$ be an
irreducible adjoint preserving $n$ dimensional representation, such
that $\Psi(\Aeps)=\varepsilon\Idn$, $\varepsilon\in\Real_+$ and
$\Psi(\Au)$ has no multiple eigenvalues. Let $\alpha$ be given by
$\tan(\tfrac12\alpha)=\varepsilon$, $0<\alpha<\pi/2$. Then there
exists $\beta\in\Real$, such that we can label the orthonormal bases
for $\Hil$ which are the eigenspaces $\Psi(\Au)$
%[
\begin{align}
\ket(\beta)\,,\,
\ket(\beta+\alpha)\,,\,
\ket(\beta+2\alpha)\,,\,
\ldots\,,\,
\ket({\beta+(n-1)\alpha}) 
\label{rep_kets}
\end{align}
%]
where
%[
\begin{align}
\Au\ket(\beta+m\alpha) &= e^{i(\beta+m\alpha)}\ket(\beta+m\alpha)
\sp
m = 0\sdotsc n-1
\label{rep_u}
\end{align}
%]
There also exists a set of complex constants
$C_{\beta+m\alpha}\in\Cmpx$ for $m=0\sdotsc n$ satisfying
%[
\begin{align}
|C_{\beta+m\alpha}|^2 &= 
\sec(\tfrac12\alpha)\cos(\beta-\tfrac12\alpha+m\alpha) + R
\label{rep_C}
\end{align}
%]
and $\CC(\beta+m\alpha)\ne0$ for $m=1\sdotsc n-1$ so that
%[
\begin{equation}
\begin{aligned}
\Aa_-\ket(\beta+m\alpha) &= \CC(\beta+m\alpha)\ket({\beta+(m-1)\alpha})
\sp &&
m = 1\sdotsc n-1
\\
\Aa_+\ket(\beta+m\alpha) &= 
\cnj{\CC({\beta+(m+1)\alpha})}\ket({\beta+(m+1)\alpha})
\sp &&
m = 0\sdotsc n-2
\end{aligned}
\label{rep_a_pm}
\end{equation}
%]

The action of $\Aa_-$ and $\Aa_+$ on the first and last vectors of
(\ref{rep_kets}) respectively are given by either
%[
\begin{align}
\Aa_-\ket(\beta) = 0 \sp 
\Aa_+\ket({\beta+(n-1)\alpha})=0
\label{rep_S2_cond}
\end{align}
%]
or
%[
\begin{align}
\Aa_-\ket(\beta) =
\CC(\beta)
\ket({\beta+(n-1)\alpha})
\sp
\Aa_+\ket({\beta+(n-1)\alpha}) =
\cnj{\CC(\beta)}
\ket(\beta)
\label{rep_T2_cond}
\end{align}
%]
In the first case $\CC(\beta)=\CC({\beta+n\alpha})=0$ satisfies
(\ref{rep_C}).
In the second case $\CC(\beta)=\CC({\beta+n\alpha})\ne0$ satisfies
(\ref{rep_C}) and
%[
\begin{align}
n\alpha=2\pi k
\label{rep_T2_cond_nk}
\end{align}
%]
where $k\in\Natn$, $1\le k<n/2$ and $n$ and $k$ are relatively prime.
\end{theorem}
%

%%%%%%%%%%%%%%%%%%%%%%%%%%%%%%%%%%%%%%%%%%%%%%%%%%%%
\begin{proof}

Let $\ket(\theta)$ be a normalised eigenvector of $\Psi(\Au)$ with
eigenvalue $\lambda$ then
$|\lambda|^2=\bra(\theta)\Au^\dagger\Au\ket(\theta)=1$ hence we can
set $\lambda=e^{i\theta}$. From (\ref{alg_def_u}) we have
$\Aw\ket(\theta)=\cos(\theta)\ket(\theta)$ and
$\Az\ket(\theta)=\sin(\theta)\ket(\theta)$. 
Let $N:\Real\mapsto\Real$ be given by
$N(\theta)=\sec(\tfrac12\alpha)\cos(\theta-\tfrac12\alpha) + R$.
From (\ref{alg_ap_am}) we have
%[
\begin{align*}
\|\Aa_-\ket(\theta)\|^2 &=
\bra(\theta)\Aa_+\Aa_-\ket(\theta)=
\bra(\theta)(\Aw+\tan(\tfrac12\Aalpha)\Az+R)\ket(\theta)=
\cos(\theta)+\tan(\tfrac12\alpha)\sin(\theta)+R = N(\theta)
\\
\|\Aa_+\ket(\theta)\|^2 &=
\bra(\theta)\Aa_-\Aa_+\ket(\theta)=
\bra(\theta)(\Aw-\tan(\tfrac12\Aalpha)\Az+R)\ket(\theta)=
\cos(\theta)-\tan(\tfrac12\alpha)\sin(\theta)+R = N(\theta+\alpha)
\end{align*}
%]
Thus if $\ket(\theta)$ is a eigenvector then $N(\theta)\ge0$ and
$N(\theta+\alpha)\ge0$. Furthermore 
$N(\theta)=0\ \Leftrightarrow\ \Aa_-\ket(\theta)=0$ and  
$N(\theta+\alpha)=0\ \Leftrightarrow\ \Aa_+\ket(\theta)=0$.  

From (\ref{alg_normal_u}) we have
%[
\begin{align*}
\Au\Aa_-\ket(\theta)=
\Aa_-\Au e^{-i\Aalpha}\ket(\theta)=
\Aa_- e^{-i\alpha}e^{i\theta}\ket(\theta)=
e^{i(\theta-\alpha)}\Aa_-\ket(\theta)
\end{align*}
%]
Hence if $N({\theta})\ne0$ then there must exist another normalised
eigenvector $\ket(\theta')$ with eigenvalue $e^{i(\theta-\alpha)}$.
If we let $\Aa_-\ket(\theta)=\CC(\theta)\ket(\theta')$ then
$|\CC(\theta)|^2=N({\theta})$. Thus $\CC(\theta)$
is determined from $N({\theta})$ up to a phase.

Similarly $\Au\Aa_+\ket(\theta')=
e^{i(\theta'+\alpha)}\Aa_+\ket(\theta')=
e^{i\theta}\Aa_+\ket(\theta')$. Since the eigenspaces of $\Psi(\Au)$
are all one dimensional then $\Aa_+\ket(\theta')$ is parallel to
$\ket(\theta)$. Thus setting $\Aa_+\ket(\theta')=D_\theta\ket(\theta)$ then 
%[
\begin{align*}
D_\theta = 
\bra(\theta')\Aa_+\ket(\theta) = 
\cnj{\bra(\theta)\Aa_-\ket(\theta')} = 
\cnj{\CC({\theta'})}
\end{align*}
%]
hence $\Aa_+\ket(\theta')=\cnj{\CC({\theta'})}\ket(\theta)$.

\vspace{1em}

{\bf Claim:} We now make the following claim.  Let $r\le n$ and there
exists a set of independent normalised eigenvectors
$\Set{\ket(\beta+m\alpha)\,|\,m=0\sdotsc r-1}$ such that
$N(\beta+m\alpha)>0$ for $m=1\sdotsc r-1$.
If either 
%[
\begin{align*}
N(\beta)>0
\textup{ and }
\Aa_-\ket(\beta)=\CC({\beta})\ket({\beta+(r-1)\alpha}) 
\textup{ for some choice of }\CC({\beta})
\end{align*}
%]
or 
%[
\begin{align*}
N(\beta)=0\,,
\textup{ and }
N(\beta+r\alpha)=0
\end{align*}
%]
then $r=n$

\vspace{0.5em}

{\bf Proof of claim:} 
Let $V_r=\textup{span}\Set{\ket(\beta+m\alpha)\,|\,m=0\sdotsc r-1}$.
Since $N(\beta+m\alpha)>0$ for $m=1\sdotsc n-1$
then we can define $\CC(\beta+m\alpha)\ne0$ for $m=1\sdotsc n-1$ 
so that $\Aa_-\ket(\beta+m\alpha) =
\CC(\beta+m\alpha)\ket({\beta+(m-1)\alpha})$ for 
$m = 1\sdotsc r-1$ and 
$\Aa_+\ket(\beta+m\alpha) = 
\cnj{\CC({\beta+(m+1)\alpha})}\ket({\beta+(m+1)\alpha})$
for $m = 0\sdotsc n-2$

For the first option 
$\Aa_-\ket(\beta)=\CC({\beta})\ket({\beta+(r-1)\alpha})$  
for so $\CC({\beta})\ne0$,
so from the argument above
$\Aa_+\ket({\beta+(r-1)\alpha})=\cnj{\CC({\beta})}\ket(\beta)$.

For the second option
$\Aa_-\ket(\beta)=0$ and
$\Aa_+\ket({\beta+(r-1)\alpha})=0$.

Hence, for both cases, the action of $\Au,\Aa_+,\Aa_-$ on $V_r$
remains in $V_r$.  Thus if $r<n$ then it is obvious that $\Psi$ can be
reduced to $V_r$.  This contradicts the irreducibility of $\Psi$ or the
dimension $\Hil$ is $n$. Hence $r=n$.

\vspace{1em}
{\bf Continuation of Proof:}

Since $N({\theta})\ge0$ for all eigenvectors $\ket(\theta)$ we have
two possibilities. Either $N({\theta})>0$ for eigenvectors
$\ket(\theta)$ or there exists a
$\ket(\beta)$ such that $N(\beta)=0$. 

Taking the first case, that $N({\theta})|>0$ all eigenvectors
$\ket(\theta)$. Choose $\beta$ so that $\ket(\beta)$ is any normalised
eigenvector of $\Psi(\Au)$. We choose the phases of $\CC(\beta+m\alpha)$,
and define $\ket(\beta+m\alpha)=\cnj{\CC({\beta+m\alpha})}^{-1}
\Aa_+\ket({\beta+(m-1)\alpha})$ for $m=1\sdotsc n-1$.

From the claim above the $\Set{\ket(\beta+m\alpha)|m=0\sdotsc n-1}$
are independent. Since $N(\beta)>0$ then $\Aa_-\ket(\beta)$ must be
parallel to one of the this set. But from the claim the only
possibility is
$\Aa_-\ket(\beta)=\CC({\beta})\ket({\beta+(n-1)\alpha})$ for some
choice of phase of $\CC({\beta})$.
i.e.  (\ref{rep_T2_cond}).

Given condition (\ref{rep_T2_cond}) then
%[
\begin{align*}
e^{i\beta}\cnj{\CC(\beta)}\ket(\beta) &=
\Au \cnj{\CC(\beta)}\ket(\beta) = 
\Au\Aa_+\ket({\beta+(n-1)\alpha})=
e^{i\alpha}\Aa_+\Au\ket({\beta+(n-1)\alpha}) \\ &=
e^{i\alpha}e^{i(\beta+(n-1)\alpha)}\Aa_+\ket({\beta+(n-1)\alpha}) = 
e^{i(\beta+n\alpha)}\cnj{\CC(\beta)}\ket(\beta)
\end{align*}
%]
Hence $e^{i n\alpha}=1$ so $n\alpha=2\pi k$ for some integer
$k$. Clearly $k\le1$ and $k<n/2$ so that $0<\alpha<\pi$. Also $k$ and
$n$ must be relatively prime so that $\Psi(\Au)$ has distinct
eigenvalues 
i.e. (\ref{rep_T2_cond_nk}).

Now consider the second possibility. That there exist a $\ket(\beta)$
such that $N(\beta)=0$ and hence $\Aa_-\ket(\beta)=0$. 
Again we choose the phases of $\CC({\beta+m\alpha})$ and define 
$\ket(\beta+m\alpha)=\cnj{\CC({\beta+m\alpha})}^{-1}
\Aa_+\ket({\beta+(m-1)\alpha})$
for $m=1\sdotsc n-1$. These are all independent and none of the
$\CC({\beta+m\alpha})=0$ by the claim above.

If $\Aa_+\ket({\beta+(n-1)\alpha})\ne0$ then by the claim above
$\Aa_+\ket({\beta+(n-1)\alpha})=\cnj{\CC({\beta+n\alpha})}\ket(\beta)$.
Hence
%[
\begin{align*}
|\CC({\beta+n\alpha})|^2=
\bra({\beta+(n-1)\alpha})\Aa_-\Aa_+\ket({\beta+(n-1)\alpha})=
\cnj{\CC({\beta+n\alpha})}\,
\bra({\beta+(n-1)\alpha})\Aa_-\ket(\beta)=
0
\end{align*}
%]
hence $\CC({\beta+n\alpha})=0$ which contradicts
$\Aa_+\ket({\beta+(n-1)\alpha})\ne0$. Hence (\ref{rep_S2_cond})
\end{proof}
%%%%%%%%%%%%%%%%%%%%%%%%%%%%%%%%%%%%%%%%

By analogy with the representations of the noncommutative sphere these
we shall call representations which satisfy (\ref{rep_S2_cond})
\defn{$S^2$-type representations} and we say a 
$S^2$-type representation is
\defn{minimal} if
%[
\begin{align}
n\alpha < 2\pi
\label{rep_S2_min}
\end{align}
%]
By analogy to the representations of the noncommutative torus we shall
call representations which satisfy (\ref{rep_T2_cond})
\defn{$T^2$-type representations} and we say a $T^2$-type representation
is \defn{minimal} if
%[
\begin{align}
n\alpha=2\pi
\label{rep_T2_min}
\end{align}
%]
We now wish to find what values of $(R,n,\alpha,\beta)$
give rise to finite dimensional irreducible representations of
$\Alg(R)$, with no multiple eigenvalues of $\Psi(\Au)$. 
A summary is given in table \ref{fig_concl} below.

%%%%%%%%%%%%%%%%%%%%%%%%%%%%%%%%%%%%%%%%%%%%%%%%%%%%%%%%%%%%%%%

\subsection{$S^2$-type representations}

Before we look at the different types of $S^2$-type representation and
when they exist we shall give some basic facts about $S^2$-type
representation.

\begin{lemma}
\label{lm_S2_basic}
If $\Psi$ is an $S^2$ representation with $(R,n,\alpha,\beta)$ given in
theorem \ref{th_reps} then we can find a basis
(\ref{rep_kets}) such that $\CC(\beta+m\alpha)\in\Real_+$.
Also there exist an equivalent representation $\tilde\Psi$ with that same
values of $(R,n,\alpha)$ and with $\beta$ replaced by $\tilde\beta$ 
where $-2\pi<\tilde\beta-\tfrac12\alpha\le0$.

If $\Psi_1$ and $\Psi_2$ are $S^2$ representation with the same
$(R,n,\alpha,\beta)$ then they are equivalent.
\end{lemma}

\begin{proof}
We can see that given $\alpha$ and $\beta$ the basis (\ref{rep_kets})
is only defined up to a phase, likewise $\CC(\beta+m\alpha)$ is also
only defined up to a phase. Thus given the set
$\Set{\nu_m\in\Cmpx\,|\,\nu_m\cnj{\nu_m}=1\,,\,m=0\sdotsc n-1}$, we
can always make the following replacement
%[
\begin{align}
\ket({\beta+m\alpha})\to 
\ket({\beta+m\alpha})' = 
\nu_m\ket({\beta+m\alpha}) 
\sp 
\CC(\beta+m\alpha)\to 
\CC(\beta+m\alpha)' = 
\nu_{m-1}\cnj{\nu_m}\CC(\beta+m\alpha)
\label{rep_phase}
\end{align}
%]
without changing the equations (\ref{rep_u}) and
(\ref{rep_a_pm}). Thus if we set
$\nu_m=\nu_{m-1}\CC(\beta+m\alpha)/|\CC(\beta+m\alpha)|$ then all the
$\CC(\beta+m\alpha)'$ are real and positive. 

Clearly replacing $\beta$ with $\tilde\beta\to\beta+2\pi k$ for 
$k\in\Natn$ doesn't change the
representation. Therefore we can place $\tilde\beta$ in any $2\pi$
range. The one chosen makes the calculations below simpler.

If $\Psi_1$ and $\Psi_2$ are representation with the same
$(R,n,\alpha,\beta)$ then, setting the $\CC(\beta+m\alpha)$ to be real
and positive, the bases (\ref{rep_kets}) are the same (up
to an overall choice phase), and the action (\ref{rep_u}),
(\ref{rep_a_pm}) on these basis elements
are the same, therefore the representations are equivalent.
\end{proof}

For the rest of this subsection we assume we are given $R$ and $n$,
and we wish to find $\alpha$ and $\beta$ so that
$\Psi(R,n,\alpha,\beta)$ is $S^2$-type irreducible representation.  
Given $\alpha$ and $\beta$ we write $\beta'=\beta-\tfrac12\alpha$.
For a general $S^2$ representation (\ref{rep_C}) and 
(\ref{rep_S2_cond}) imply
%[
\begin{align}
\cos(\beta')=\cos(\beta'+n\alpha)=-R\cos(\tfrac12\alpha)
\label{rep_lm_S2_cond}
\end{align}
%]
The first equality is solved by setting $\beta'+n\alpha=2\pi k \pm
\beta'$ for some $k\in\Natn$. This implies either 
%[
\begin{align}
\alpha=2\pi k/n
\label{rep_lm_S2_cond1}
\end{align}
%]
or 
%[
\begin{align}
\beta'=\pi k -\tfrac12 n \alpha
\label{rep_lm_S2_cond2}
\end{align}
%]
We can now place some simple constraints on $(R,n,\alpha,\beta)$ such
that $\Psi(R,n,\alpha,\beta)$ is an irreducible $S^2$ representation.
From (\ref{rep_lm_S2_cond}) we see that for $\Psi$ to be an $S^2$
representation then  
%[
\begin{align}
R\le\sec(\tfrac12\alpha)=(1+\varepsilon^2)^{1/2}
\end{align}
%]
since $\varepsilon=\tan(\tfrac12\alpha)$. By looking at
(\ref{rep_C}) we see that $\Psi(R,n,\alpha,\beta)$ is an irreducible 
$S^2$ representation if and only if
%[
\begin{align}
\cos(\beta'+m\alpha) + R\cos(\tfrac12\alpha) > 0
\sp
m=1\sdotsc n-1
\label{rep_S2_ineq}
\end{align}

%%%%%%%%%%%%%%%%%%%%%%%%%%%%%%%%%%%%%%%%%%%%%%%%%%%%]
\subsubsection*{Minimal $S^2$ representation}

The first result is on the existence and uniqueness of minimal
$S^2$-type representations.
\begin{lemma}
\label{lm_S2min}
Given $R$\, and $n$ there exists minimal $S^2$-type representation if and only
if \mbox{$-1<R<\sec(\pi/n)$}. This representation is unique and is given by
%[
\begin{align}
\cos(\tfrac12n\alpha)+R\cos(\tfrac12\alpha)=0
\sp
0<\alpha<2\pi/n
\sp
\beta'=-\tfrac12 n\alpha
\label{rep_minS2}
\end{align}
%]
\end{lemma}

\begin{proof}
If $\Psi$ is minimal we must exclude (\ref{rep_lm_S2_cond1}).
Since we choose $\beta'$ so that $-2\pi<\beta'\le0$ from
(\ref{rep_lm_S2_cond2}) we have $-2<k-n\alpha/(2\pi)\le0$. Applying
$0<n\alpha/(2\pi)<1$ we have
%[
\begin{align*}
-2<n\alpha/(2\pi)<k\le 1+n\alpha/(2\pi) < 1
\end{align*}
%]
therefore $k=0$ or $-1$. 

We shall exclude the case $k=-1$. If $k=-1$ so that
$\beta'=-\pi-\tfrac12 n\alpha$, therefore $-2\pi<\beta'<-\pi$. Now
$\beta'+\alpha=-\pi-\tfrac12(n-2)\alpha\le-\pi$ since $n\ge2$. Thus
$-2\pi<\beta'<\beta'+\alpha\le -\pi$. Since $\cos$ is strictly
decreasing in the range $-2\pi\ldots-\pi$ we have
$-1\le\cos(\beta'+\alpha)<\cos(\beta')<1$. Thus
$\cos(\beta'+\alpha)-\cos(\beta')=\cos(\beta'+\alpha)+R\cos(\tfrac12\alpha)<0$
which contradicts (\ref{rep_S2_ineq}).

Thus if $\Psi$ is a minimal $S^2$ representation we have
(\ref{rep_minS2}). 

Now consider the function
$\hat{R}(\alpha)=-\cos(\tfrac12n\alpha)/\cos(\tfrac12\alpha)$ for the
range $0<\alpha<2\pi/n$. 
We observe that $\hat{R}(0)=-1$ and $\hat{R}(2\pi/n)=\sec(\pi/n)$ and
%[
\begin{align*}
\hat{R}'(\alpha)&=
(\cos(\tfrac12\alpha))^{-2}(
\tfrac12n\alpha\cos(\tfrac12\alpha)\sin(\tfrac12n\alpha)-
\tfrac12\alpha\cos(\tfrac12n\alpha)\sin(\tfrac12\alpha))
\\
&=
(\cos(\tfrac12\alpha))^{-2}
\tfrac12\alpha((n-1)\cos(\tfrac12\alpha)\sin(\tfrac12n\alpha)
+ \sin(\tfrac12(n-1)\alpha) ) > 0
\end{align*}
%]
Thus
$\hat{R}:\Set{\alpha|0<\alpha<2\pi/n}\mapsto\Set{R|-1<R<\sec(\pi/n)}$ 
is an invertible function so we can uniquely solve
(\ref{rep_minS2}). This also implies that $-1<R<\sec(\pi/n)$.

We now wish to show that an $\alpha,\beta'$ satisfying
(\ref{rep_minS2}) is a representation. This simply requires showing
(\ref{rep_S2_ineq}) is satisfied. Consider separately the cases
$m=1\sdotsc\lfloor\tfrac12n\rfloor$ and
$m=\lfloor\tfrac12n\rfloor+1\sdotsc n-1$. In the first case we have
$-\tfrac12n<m-\tfrac12n\le0$ this implies $-\pi\le-\tfrac12n\alpha<
(m-\tfrac12n)\alpha \le0$.  Since $\cos$ is strictly increasing in this
range we have
$-1\le\cos(-\tfrac12n\alpha)<\cos(m\alpha-\tfrac12n\alpha)\le1$ and
hence (\ref{rep_S2_ineq}) is satisfied.  Likewise if
$m=\lfloor\tfrac12n\rfloor+1\sdotsc n-1$ then $0\le
m-\tfrac12n<\tfrac12n$ so $0\le
m\alpha-\tfrac12n\alpha < \tfrac12n\alpha \le \pi$.  Since $\cos$ is
strictly decreasing in this range we have
$-1\le\cos(\tfrac12n\alpha)<\cos(m\alpha-\tfrac12n\alpha)\le1$ and
hence (\ref{rep_S2_ineq}) is satisfied.
\end{proof}

%%%%%%%%%%%%%%%%%%%%%%%%%%%%%%%%%%%%%%%%%%%%%%%%%%%%%%%

\begin{figure}
\setlength{\unitlength}{6cm}
\begin{picture}(1,1.1)(0,-0.1)
\put(0,0){
\includegraphics[height=\unitlength]{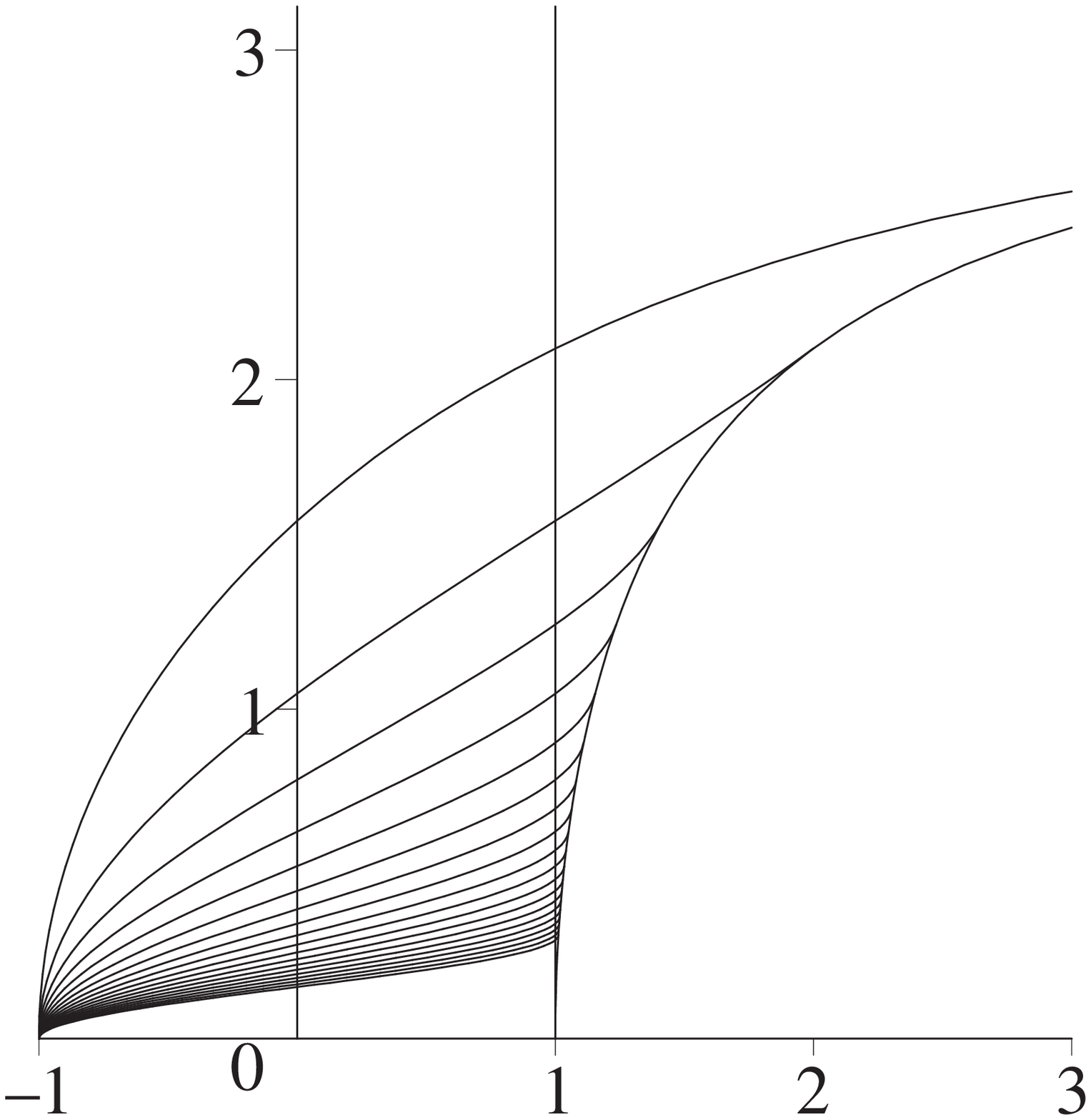}
}
\put(0.2,0.8){$\alpha$}
\put(0.8,0.0){$R$}
\put(0.7,0.6){\makebox(0,0)[lt]{$R=\sec(\alpha)$}}
\put(0,0){\makebox(0,0)[lt]{Figure \ref{fig_min_S2}a}}
\end{picture}
\qquad
\setlength{\unitlength}{5cm}
\begin{picture}(1,1.2)(0,-0.2)
\put(0,0){
\includegraphics[height=\unitlength]{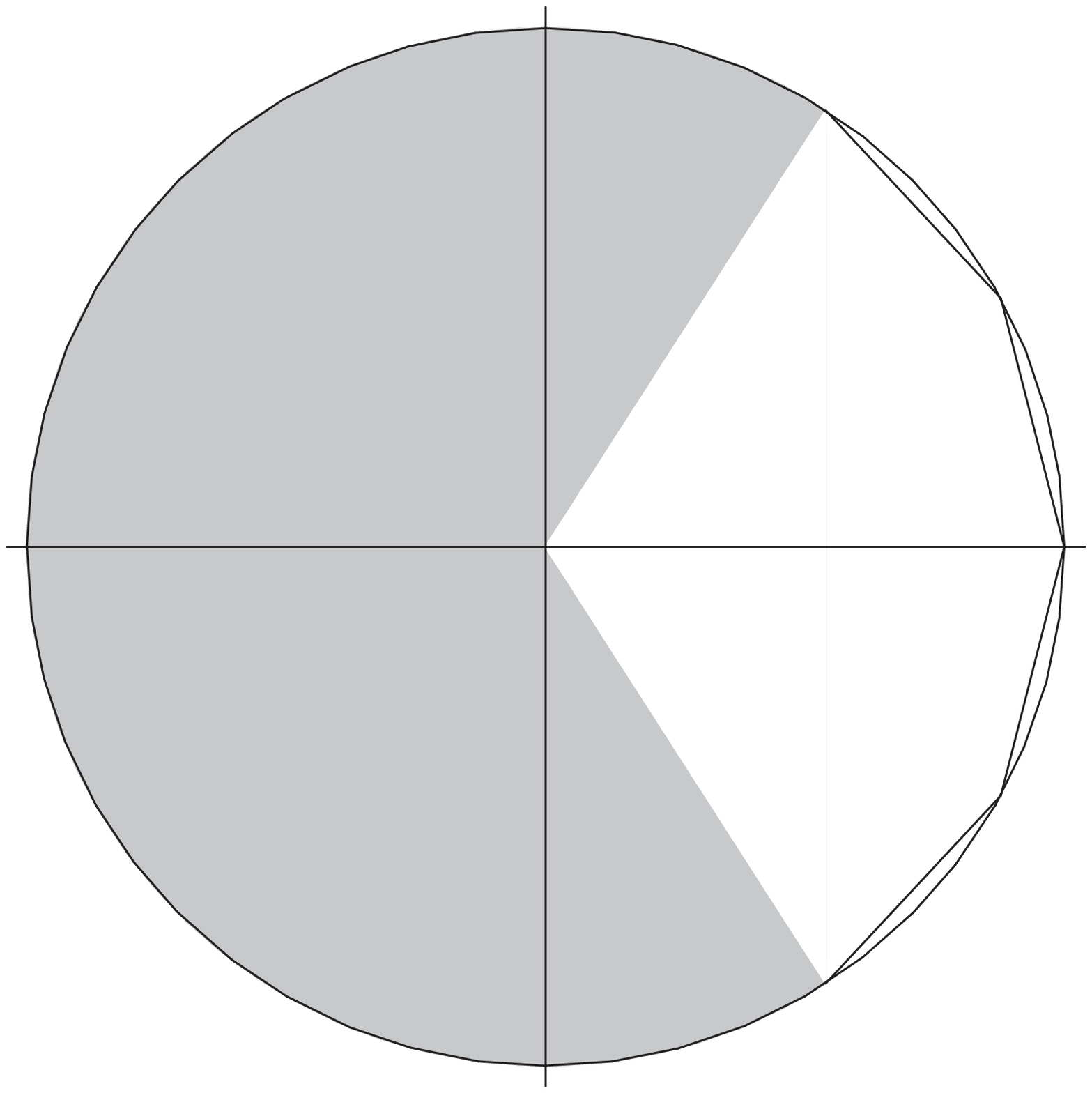}}
\put(0,0){\makebox(0,0)[tl]
{\tabelm{Figure \ref{fig_min_S2}b: $n=4,\,\alpha=.5,$
         \\$R = -.557,\,\beta=-1$}}}
\put(0.1,0.55){\tabelm{forbidden\\sector}}
\end{picture}
\qquad
\begin{picture}(1,1.2)(0,-0.2)
\put(0,0){
\includegraphics[height=\unitlength]{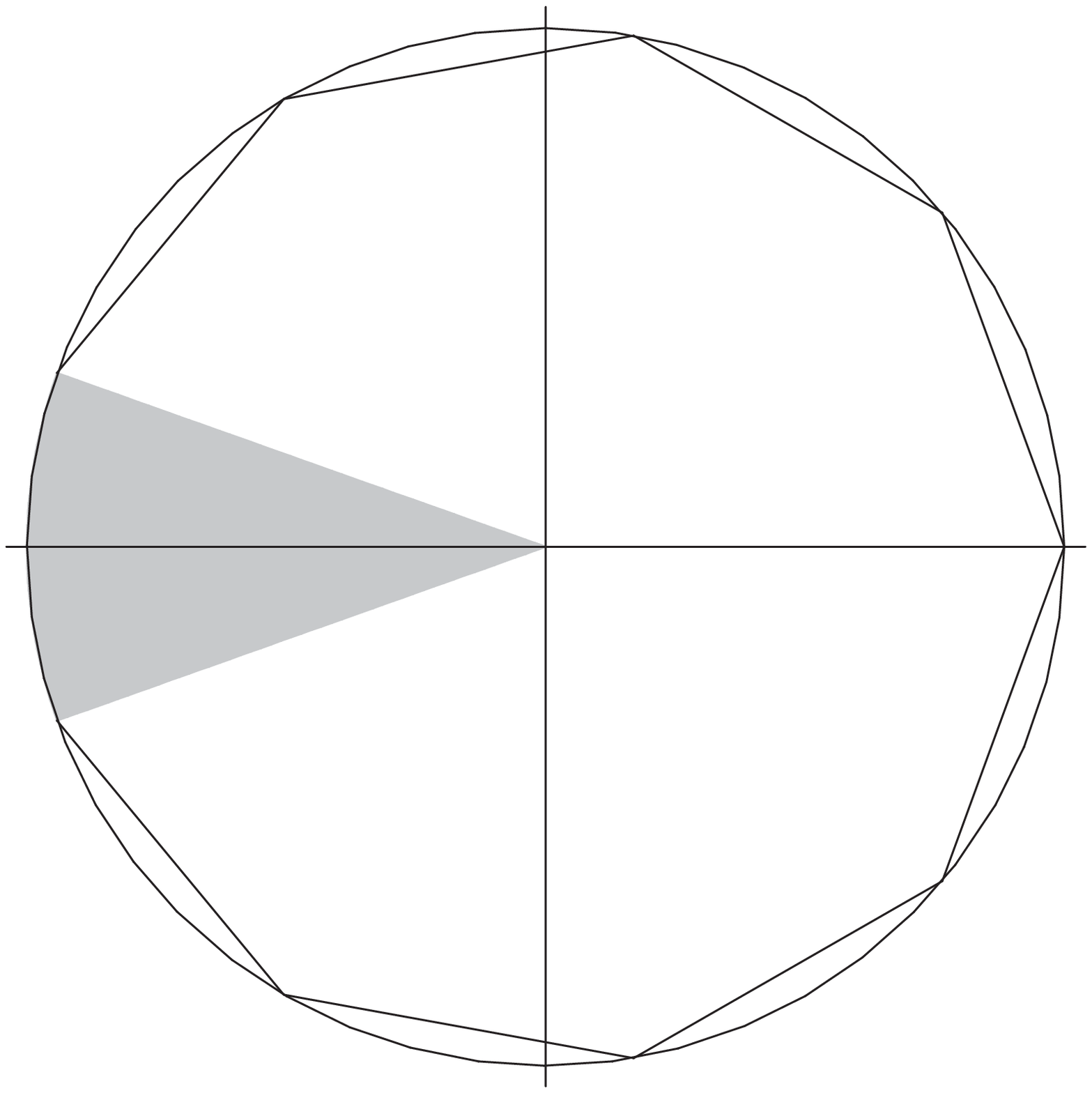}}
\put(0,0){\makebox(0,0)[tl]
{\tabelm{Figure \ref{fig_min_S2}c: $n = 8\,, \alpha = .7,$
                            \\$R = 1.003,\, \beta = -2.8$}}}
\put(0.,0.55){\tabelm{forbidden\\sector}}
\end{picture}
\caption{The minimal $S^2$ representation. Figure \ref{fig_min_S2}a
gives the relationship between $R$, $n$ and $\alpha$ for
$n=2,3,\ldots$ Figures \ref{fig_min_S2}b and \ref{fig_min_S2}c are
represent different representations.}
\label{fig_min_S2}
\end{figure}

The minimal $S^2$ representation may be said to be 
the closest to the representation
of the true noncommutative sphere, since they exist for $-1<R<1$ when
the commutative limit is a topologically the sphere, 
they are uniquely determined for each $R$
by $n$, and $\varepsilon=\tan(\tfrac12\alpha)\to0$ as $n\to\infty$. 

As lemma \ref{lm_S2min} shows there exists a minimal
$S^2$ representation of dimension $n$ for $-1<R<\sec(\pi/n)$. The
relationship between $R$, $n$ and $\alpha$ is shown in Figure
\ref{fig_min_S2}a,  $n=2,3,\ldots$. Given a permitted $R$ and $n$ we
can draw a diagram to represent $\Psi$. This an open polygon inside a
circle. The vertices are at the points
$e^{i(\beta'+m\alpha)}$. Equations (\ref{rep_S2_ineq}) implies that
there is a sector (called the ``forbidden sector'') and all the
vertices must lie to the right of this sector. The open polygon begins
and ends on the edge of this sector. Figures \ref{fig_min_S2}b and
\ref{fig_min_S2}c are represent representations for different $n$ and
$R$. Note for Figure \ref{fig_min_S2}c, $R>1$. 

%%%%%%%%%%%%%%%%%%%%%%%%%%%%%%%%%%%%%%%%%%%%%%%%%%%%%%%
\subsubsection*{Non minimal $S^2$ representation}

The situation for non minimal $S^2$ representation is more
complicated. We have the following lemmas: 

\begin{lemma}
\label{lm_S2_non_min_R}
If $\Psi$ is a non minimal $S^2$-type representation then $R\ge1$.
\end{lemma}
\begin{proof}
If $\ket(\theta)$ is a basis vector then from (\ref{rep_S2_ineq})
$\cos(\theta)\ge-R\cos(\tfrac12\alpha)$. Therefore if we set
$\theta'=(\theta)\!\!\!\!\mod 2\pi$ so that $0\le\theta'<2\pi$, then $\theta'$
must lie either in the range $0\le\theta'\le\pi-\tfrac12\delta$ or
$\pi+\tfrac12\delta\le\theta'<2\pi$ where
$\cos(\tfrac12\delta)=R\cos(\tfrac12\alpha)$ and $0<\delta$. 
The range
$\Set{\theta'|\pi-\tfrac12\delta<\theta'<\pi+\tfrac12\delta}$ is
called the forbidden sector.

Since $\Psi$ is non minimal there exists an $m$ such that
$0\le(\beta'+m\alpha)\!\!\!\!\mod2\pi\le\pi-\tfrac12\delta$ and
$\pi+\tfrac12\delta\le(\beta'+(m+1)\alpha)\!\!\!\!\mod 2\pi<2\pi$. Thus
$0<\delta\le\alpha<\pi$. Since $\cos$ is decreasing we have
$0<\cos(\tfrac12\alpha)\le\cos(\tfrac12\delta)<1$. Hence result.
\end{proof}

\begin{figure}
\begin{center}
\setlength{\unitlength}{6cm}
\begin{picture}(1,1.1)(0,0)
\put(0,0){
\includegraphics[height=\unitlength]{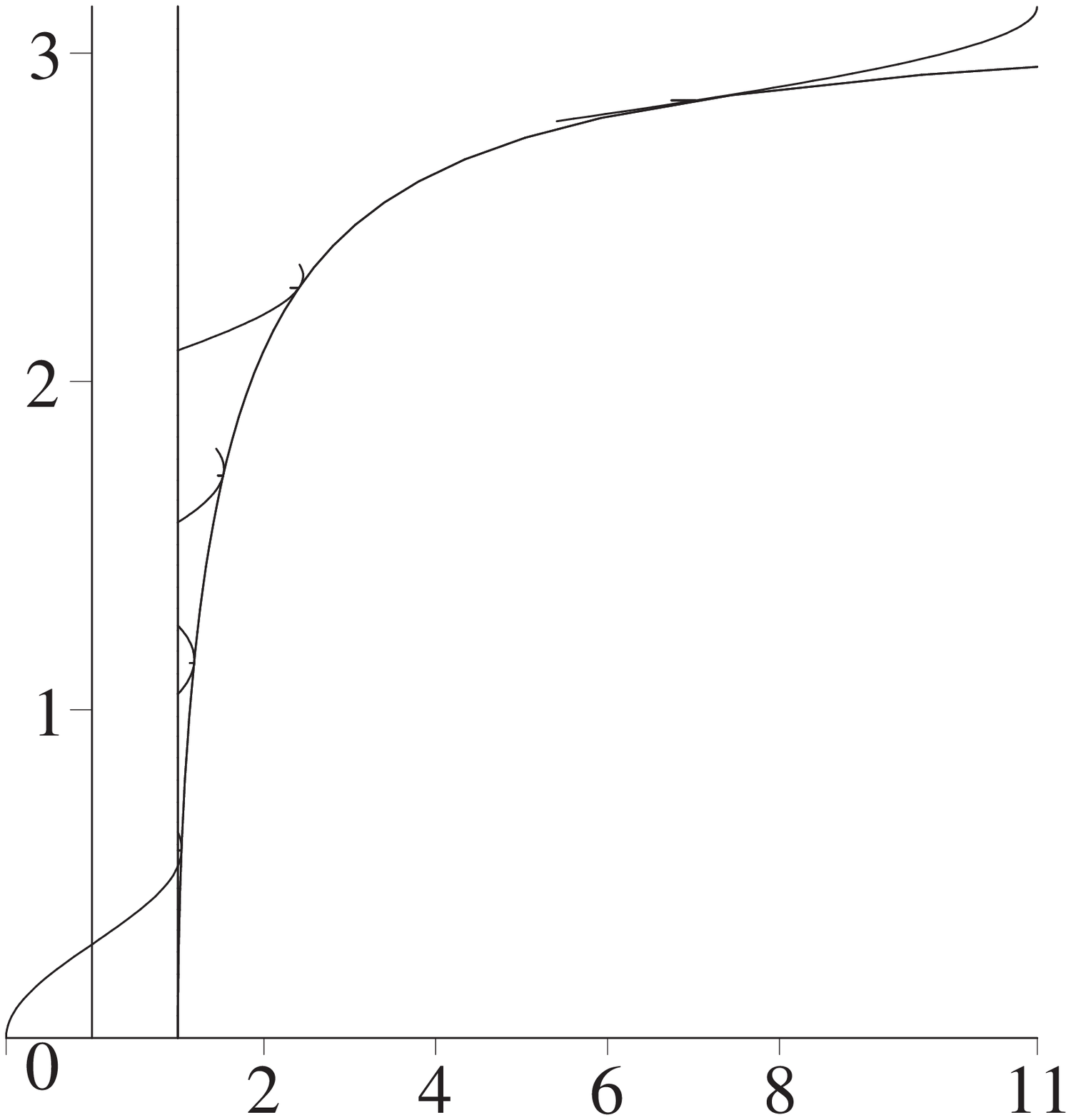}}
\put(0.0,0.8){$\alpha$}
\put(0.8,0.0){$R$}
\put(0.3,0.7){\makebox(0,0)[lt]{$R=\sec(\alpha)$}}
\put(0,0){\makebox(0,0)[lt]{Figure \ref{fig_non_min_S2}a}}
\end{picture}
\qquad
\setlength{\unitlength}{6cm}
\begin{picture}(1,1.1)(0,0)
\put(0,0){
\includegraphics[height=\unitlength]{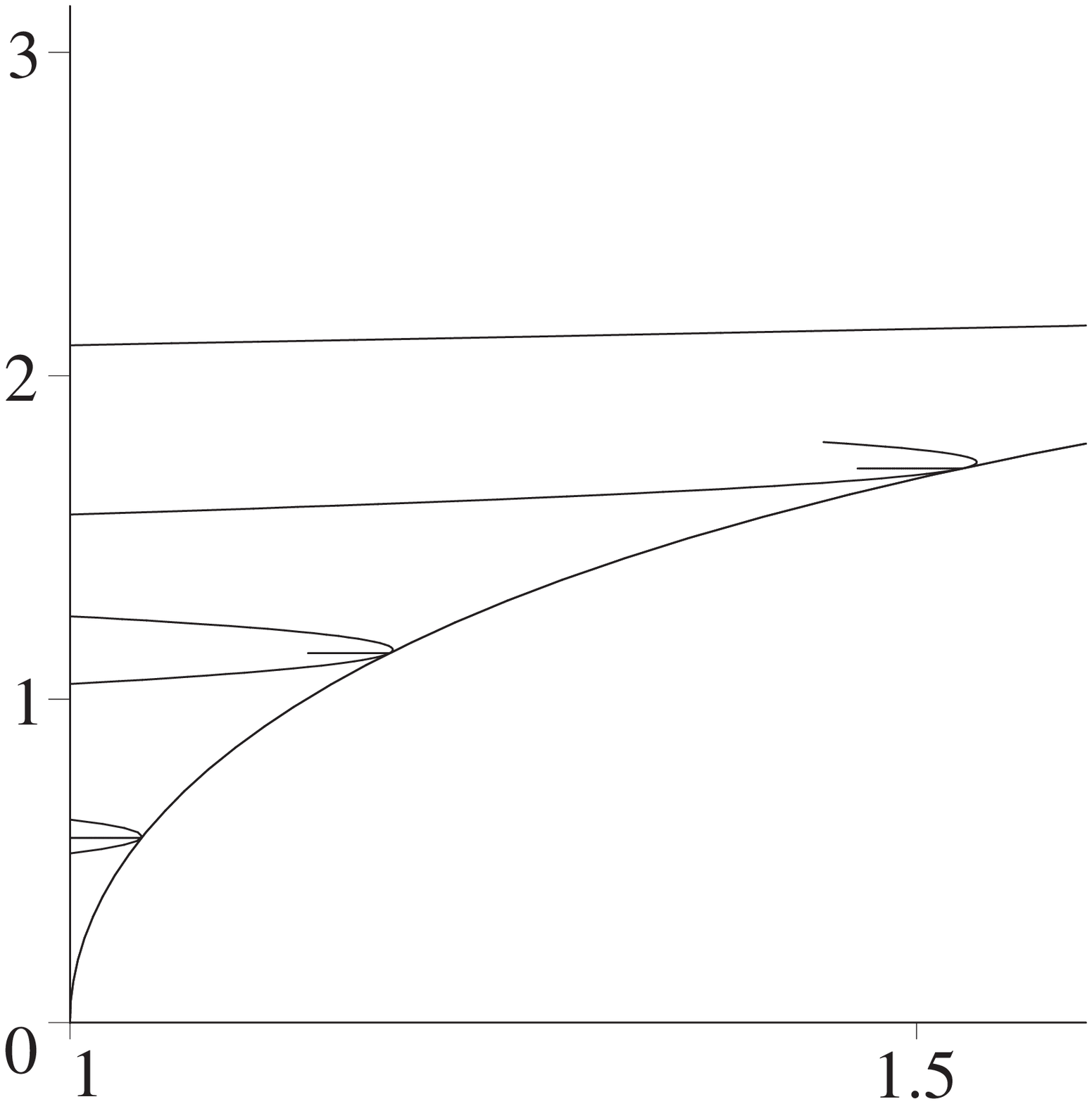}}
\put(0.0,0.8){$\alpha$}
\put(0.95,0.0){$R$}
\put(0.6,0.3){\makebox(0,0)[lt]{$R=\sec(\alpha)$}}
\put(0,0){\makebox(0,0)[lt]{Figure \ref{fig_non_min_S2}b}}
\end{picture}
\end{center}
\caption{Non minimal $S^2$ representation, the relationship between
$R$ and $\alpha$ for $n=11$}
\label{fig_non_min_S2}
\end{figure}

\begin{figure}[t]
\setlength{\unitlength}{4.3cm}
\begin{picture}(1,1.3)(0,-0.3)
\put(0,0){
\includegraphics[height=\unitlength]{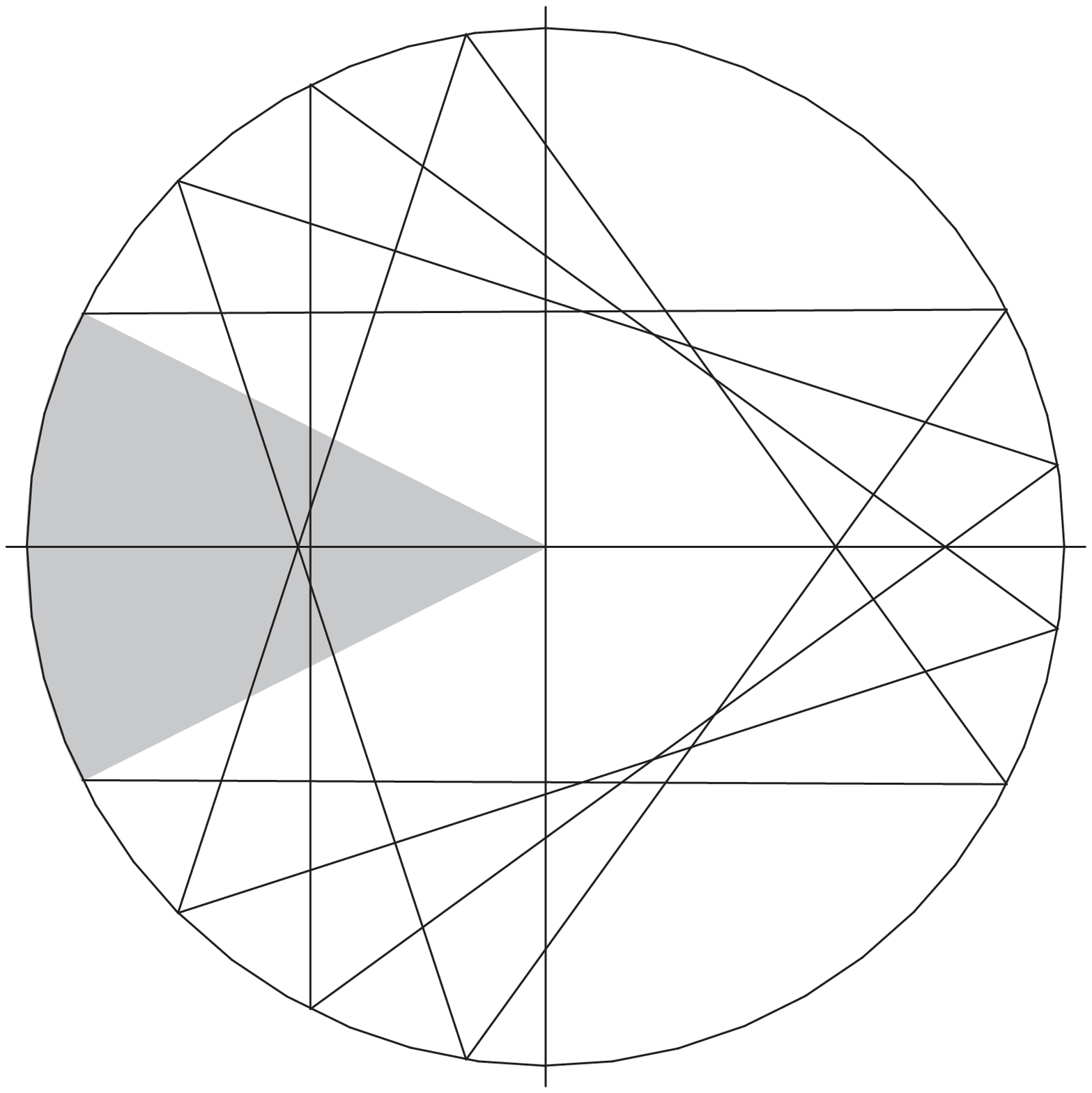}}
\put(0,0){\makebox(0,0)[tl]
{\tabelm{
Figure \ref{fig_S2_circ}a: (good) \\
$N=11,\alpha=2.20$,\\$R = 1.97,\beta=-2.67$}}}
\put(0.08,0.52){FS}
\end{picture}
\begin{picture}(1,1.3)(0,-0.3)
\put(0,0){
\includegraphics[height=\unitlength]{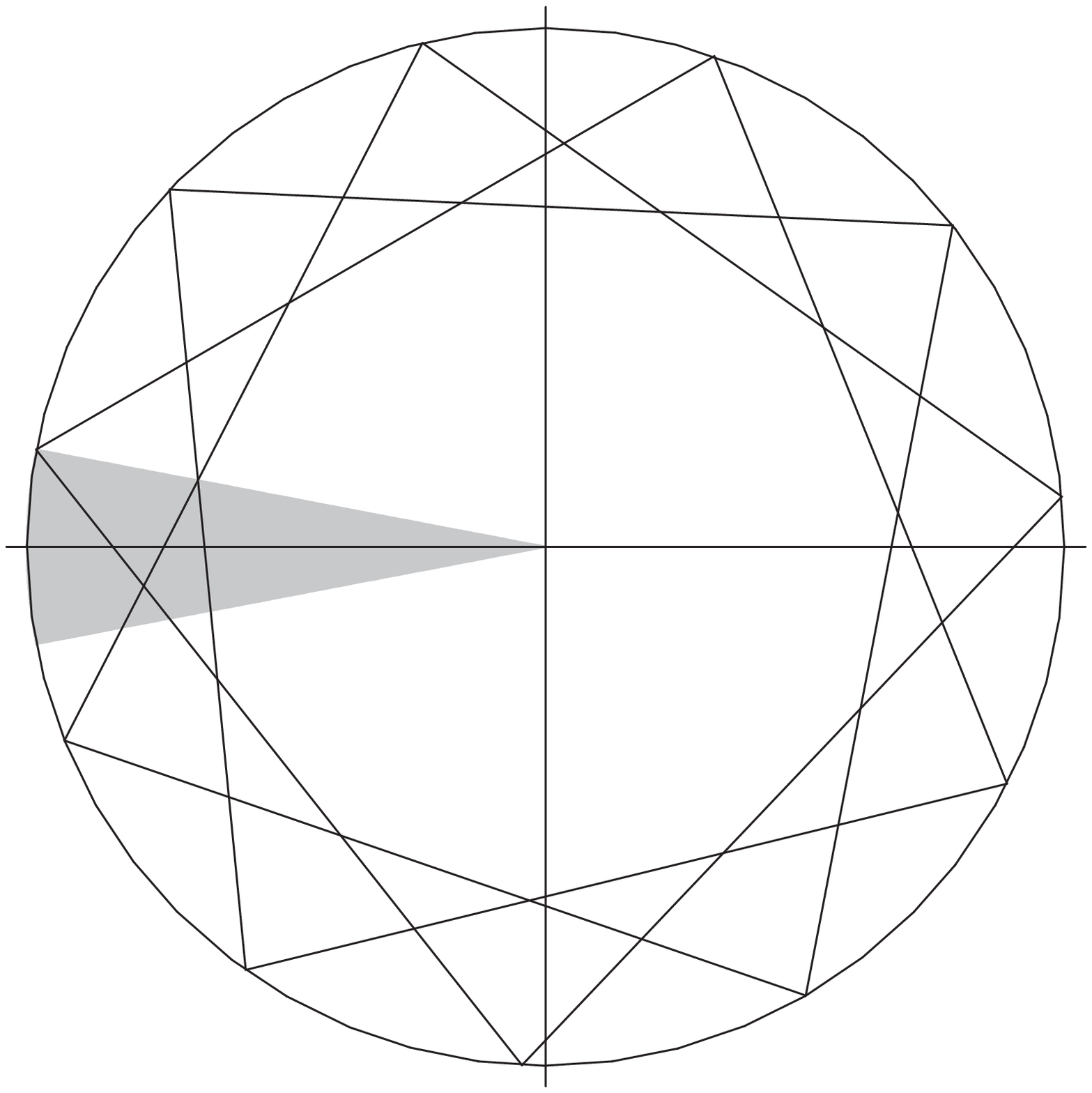}}
\put(0,0){\makebox(0,0)[tl]
{\tabelm{
Figure \ref{fig_S2_circ}b: (good) \\
$N=11,\alpha=2\pi (3/11)$,\\$R = 1.50,\beta=-2.95$}}}
\put(0.25,0.48){FS}
\end{picture}
\begin{picture}(1,1.3)(0,-0.3)
\put(0,0){
\includegraphics[height=\unitlength]{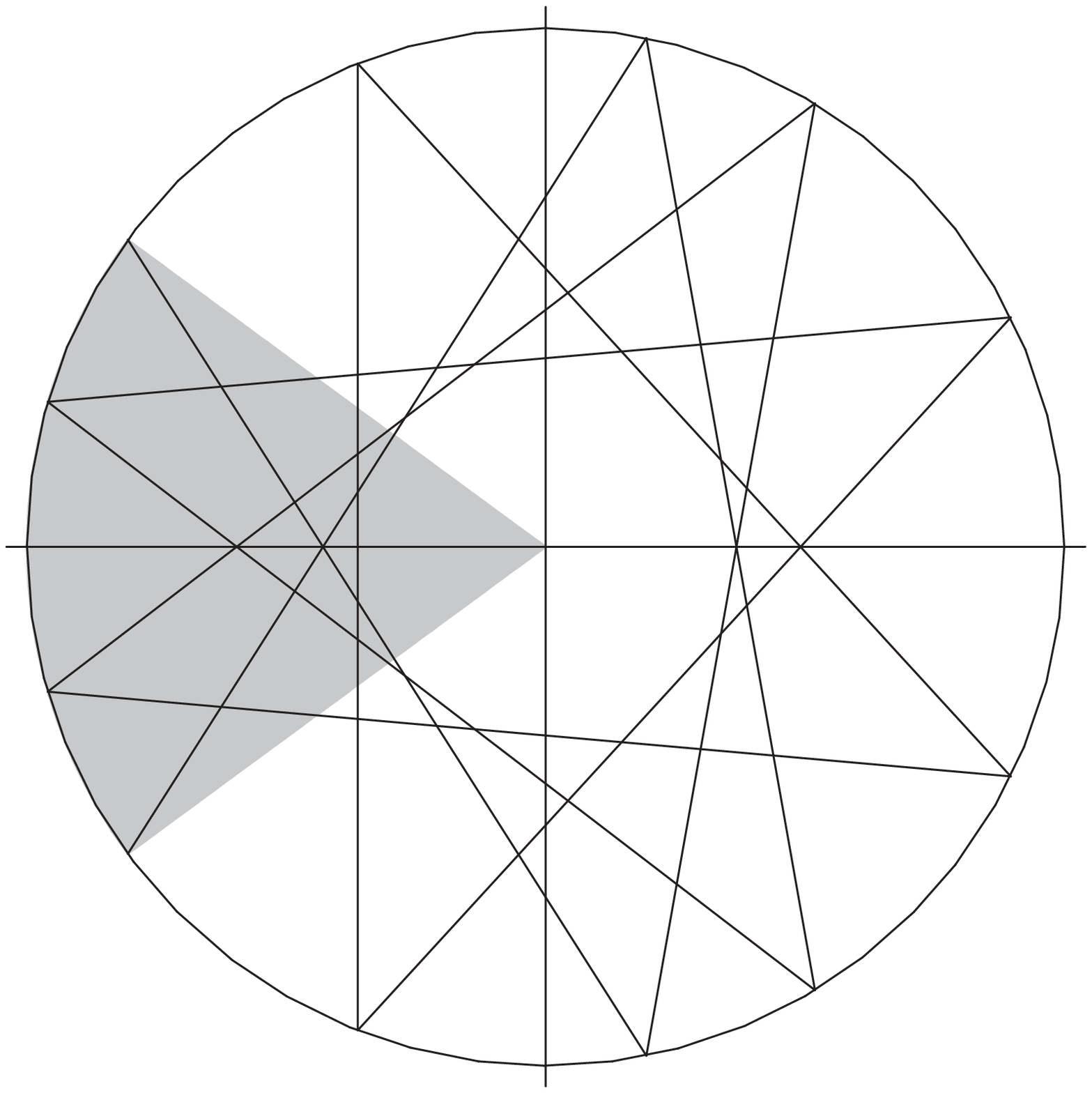}}
\put(0,0){\makebox(0,0)[tl]
{\tabelm{
Figure \ref{fig_S2_circ}c: (bad) \\
$N=11,\alpha=2.40$,\\$R = 2.22,\beta=-3.77$}}}
\put(0.05,0.48){FS}
\end{picture}
\begin{picture}(1,1.3)(0,-0.3)
\put(0,0){
\includegraphics[height=\unitlength]{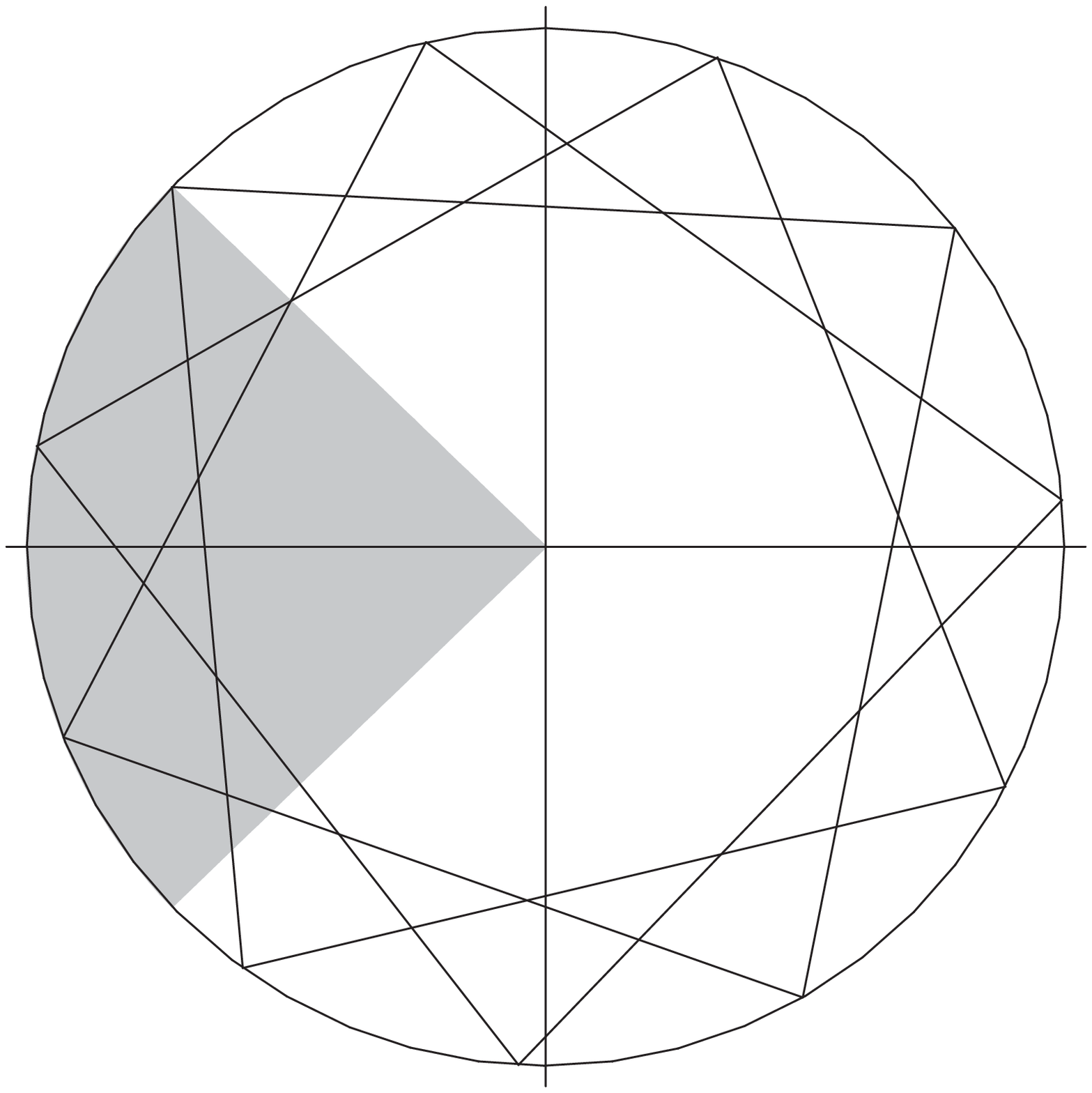}}
\put(0,0){\makebox(0,0)[tl]
{\tabelm{
Figure \ref{fig_S2_circ}d: (bad) \\
$N=11,\alpha=2\pi (3/11)$,\\$R = 1.10,\beta=-2.37$}}}
\put(0.3,0.53){FS}
\end{picture}
\caption{Permissible and non permissible values of
$n,R,\alpha,\beta$. (FS=forbidden sector)}
\label{fig_S2_circ}
\end{figure}

\begin{lemma}
If $R\le-1$ there exist no $S^2$-type representations.
\end{lemma}
\begin{proof}
Follows from lemmas \ref{lm_S2min} and \ref{lm_S2_non_min}.
\end{proof}

\begin{lemma}
\label{lm_S2_non_min}
If $R\ge1$ then there may exists many non minimal $S^2$-type
representation. Given such a $\Psi(R,n,\alpha,\beta)$ then
$1\le R\le\sec(\tfrac12\alpha)$ and
$-\tfrac32\pi<\beta'<-\tfrac12\pi$ and also
$(R,n,\alpha,\beta)$ obey either
%[
\begin{align}
&
\cos(\tfrac12n\alpha) = (-1)^{k+1} R\cos(\tfrac12\alpha)
\sp
2\pi/n<\alpha<\pi/2
\sp
\beta'=\pi k - \tfrac12 n\alpha
\sp
k = [n\alpha/(2\pi)]-1
\label{rep_S2_a}
\end{align}
%]
or
%[
\begin{align}
\alpha = 2k\pi/n
\sp
1\le k \le (n-1)/2
\sp
\cos(\beta')=-R\cos(\pi k/n)
\label{rep_S2_b}
\end{align}
%]
\end{lemma}
\begin{proof}
From (\ref{rep_lm_S2_cond}) and $|\cos(\beta')|\le1$ we have 
$R\le\sec(\tfrac12\alpha)$. Also since $R>0$ and
$\cos(\tfrac12\alpha)>0$ then $\cos(\beta')<0$ then
$-\tfrac32\pi<\beta'<-\tfrac12\pi$.

From (\ref{rep_lm_S2_cond}) we have either (\ref{rep_lm_S2_cond1}) or
(\ref{rep_lm_S2_cond2}) is true. If (\ref{rep_lm_S2_cond2}) is true
then from (\ref{rep_lm_S2_cond}) the first equation in
(\ref{rep_S2_a}) is true. Also $\beta'/\pi=k-n\alpha/(2\pi)$ so
$[\beta'/\pi]=-1=k-[n\alpha/(2\pi)]$.

If (\ref{rep_lm_S2_cond1}) is satisfied, then from
(\ref{rep_lm_S2_cond}) equation (\ref{rep_S2_b}) is satisfied.  
\end{proof}

We note that the converse is to lemma \ref{lm_S2_non_min} is not
true. That is given $(R,n,\alpha,\beta)$ which satisfies either
(\ref{rep_S2_a}) or (\ref{rep_S2_b}) there need not be a corresponding
representation $\Psi(R,n,\alpha,\beta)$. This is due to the
requirement that (\ref{rep_S2_ineq}) is satisfied.  For example given
$(R,n,\alpha,\beta)$ which is a solution to (\ref{rep_S2_a}) such that
$n$ is even and $k$ is old, then for $m=n/2$ we have
%[
\begin{align*}
\cos(\beta'+m\alpha) + R\cos(\tfrac12\alpha)=
\cos(\pi k - \tfrac12 n \alpha + m \alpha) - \cos(\beta') =
\cos(\pi k)-\cos(\beta') = -1-\cos(\beta') < 0
\end{align*}
%]
hence (\ref{rep_S2_ineq}) is not satisfied.

In general that exact set of $(R,n,\alpha,\beta)$ which have non
minimal $S^2$ representation is complicated. In figure
\ref{fig_non_min_S2} we see the relationship between $R$ and $\alpha$
for $n=11$. Figure \ref{fig_non_min_S2}b is simply a smaller
region of $R$. Figure \ref{fig_S2_circ} gives four different representation
with $n=11$, the first two are acceptable since (\ref{rep_S2_ineq}) is
satisfied, whereas the second two are unacceptable since
(\ref{rep_S2_ineq}) is not satisfied. Figures \ref{fig_S2_circ}a and
\ref{fig_S2_circ}c correspond to (\ref{rep_S2_a}). Figures
\ref{fig_S2_circ}b and \ref{fig_S2_circ}d correspond to
(\ref{rep_S2_b})

%%%%%%%%%%%%%%%%%%%%%%%%%%%%%%%%%%%%%%%%%%%%%%%%%%

\begin{figure}
\setlength{\unitlength}{5cm}
\begin{picture}(1.1,1.4)(-0.1,-0.4)
\put(0,0){
\includegraphics[height=\unitlength]{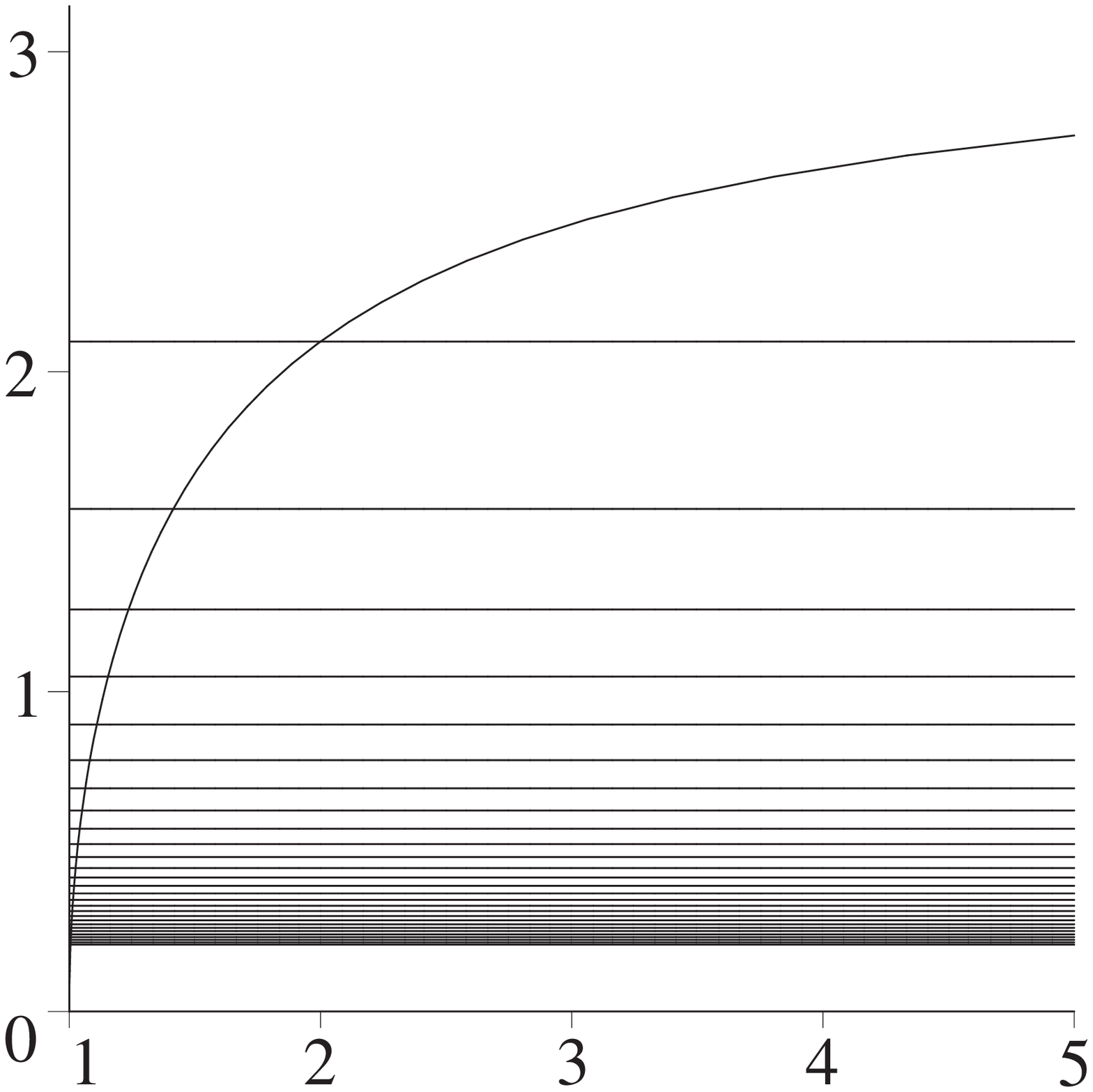}}
\put(0.0,0.8){$\alpha$}
\put(0.9,-0.01){$R$}
\put(1,1){\makebox(0,0)[rt]{$R=\sec(\alpha)$}}
\put(0.1,-0.05){\makebox(0,0)[lt]
       {\tabelm{Figure \ref{fig_T2}a: \\ 
        Minimal $T^2$ representation,\\ 
        Relationship between $R,n,\alpha$, \\
        for $n=2,3,\ldots$}}}
\end{picture}
\qquad
\begin{picture}(1.1,1.3)(-0.1,-0.4)
\put(0,0){
\includegraphics[height=\unitlength]{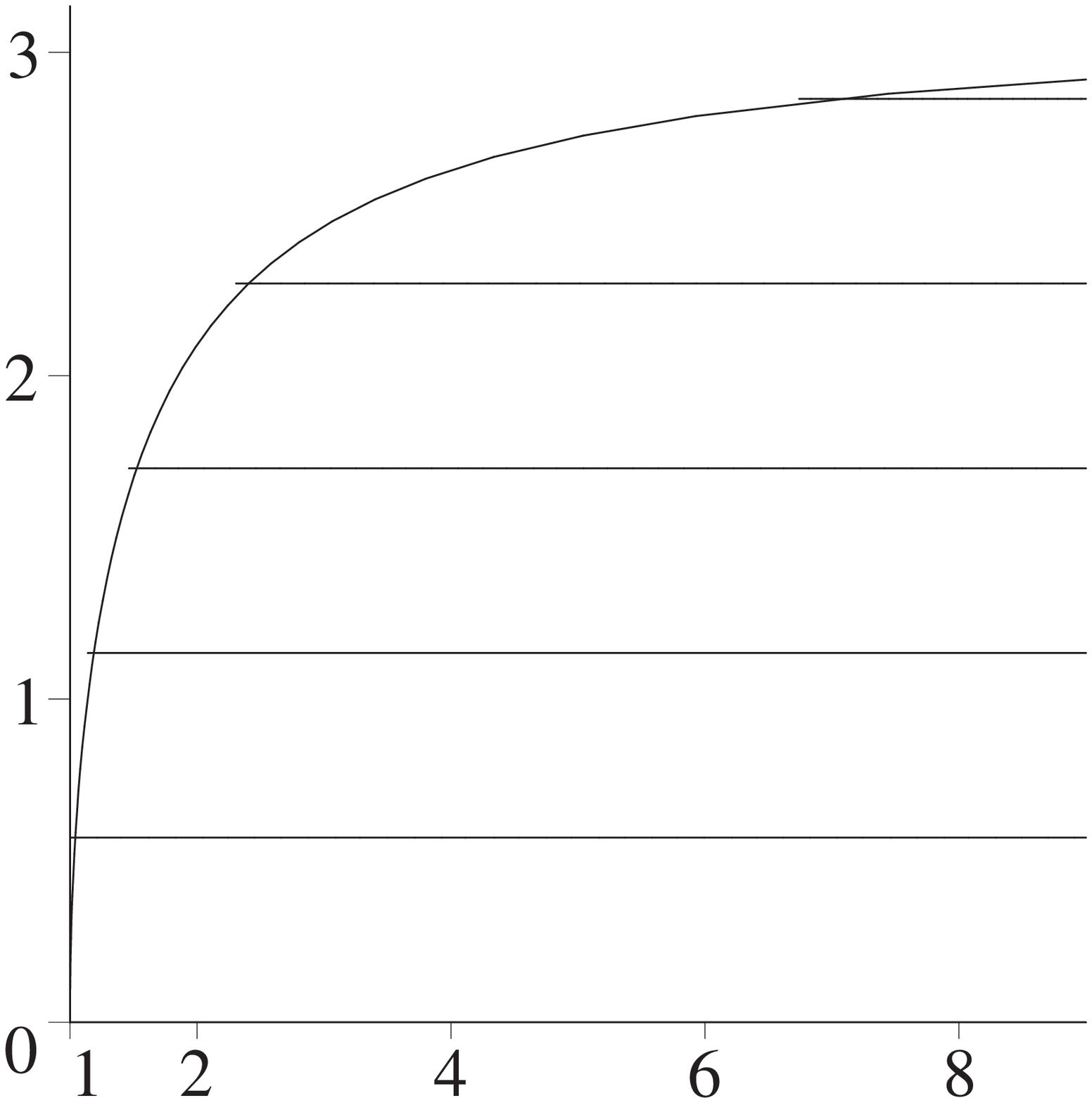}}
\put(0.0,0.8){$\alpha$}
\put(0.95,-0.01){$R$}
\put(1,1){\makebox(0,0)[rt]{$R=\sec(\alpha)$}}
\put(0.1,-0.05){\makebox(0,0)[lt]
         {\tabelm{Figure \ref{fig_T2}b: \\
         Relationship between $R,\alpha$ \\ for $n=11$.}}}
\end{picture}
\qquad
\begin{picture}(1,1.3)(0,-0.4)
\put(0,0){
\includegraphics[height=\unitlength]{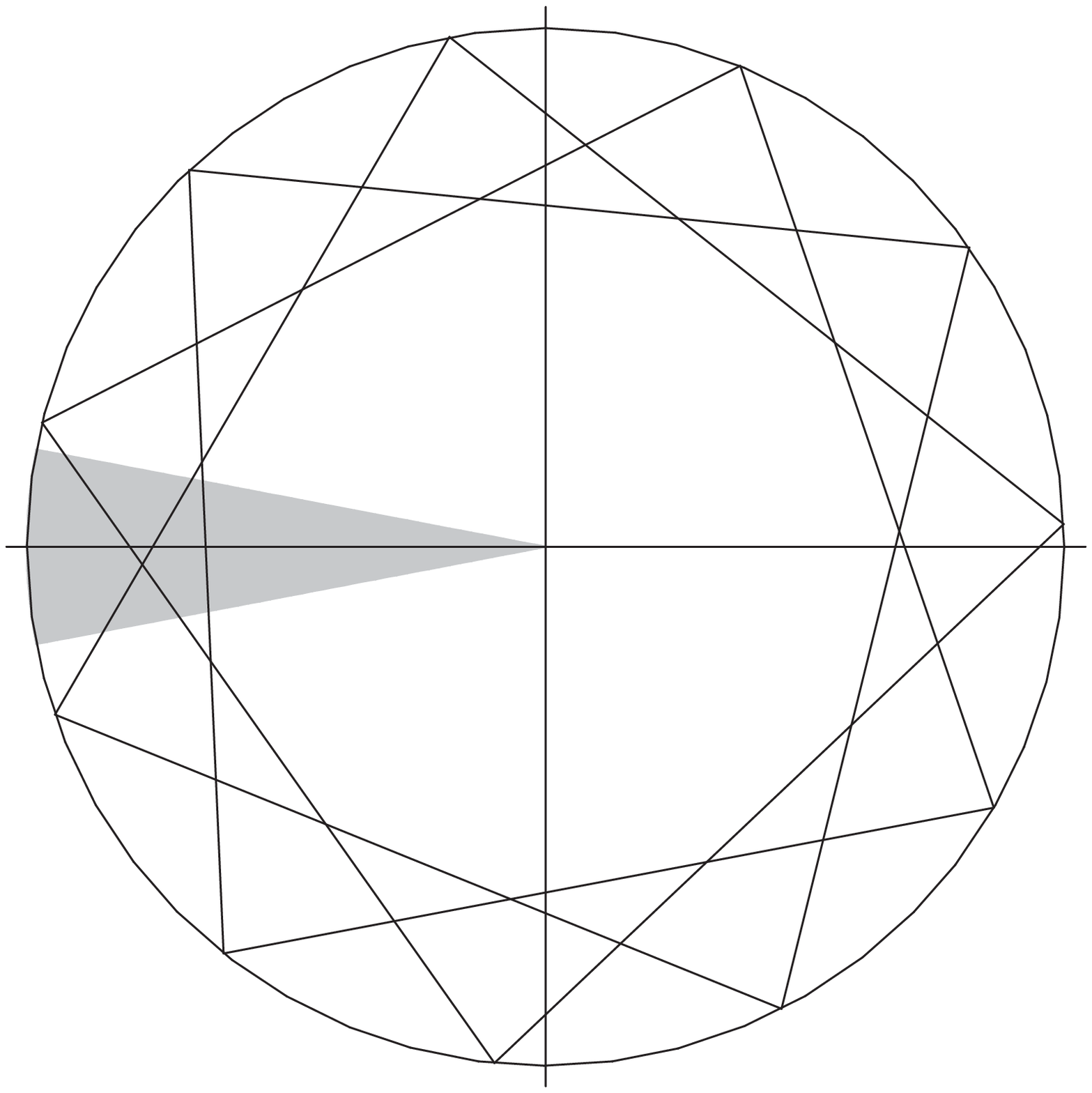}}
\put(0,0){\makebox(0,0)[tl]
{\tabelm{
Figure \ref{fig_T2}c: \\
$N=11,\alpha=2\pi (3/11)$,\\$R = 1.10,\beta=-2.9$}}}
\put(0.25,0.48){FS}
\end{picture}
\caption{Permissible region for $T^2$ representation and an example.}
\label{fig_T2}
\end{figure}

\subsection{$T^2$-type representations}

Again before we look at the different types of $T^2$-type
representation and when they exist, we shall give some basic facts
about $T^2$-type representation. Recall that for a $T^2$-type
representation $\alpha=2\pi k/n$.

We note take we can still make the replacement given by
(\ref{rep_phase}) for $m=0\sdotsc n-1$ with
$\nu_{-1}=\nu_{n-1}$. However in this case the constant 
$C_{\textup{prod}} = \prod_{m=0}^{n-1} \CC(\beta+m\alpha)$
is unchanged under the replacement
(\ref{rep_phase}). i.e. $C_{\textup{prod}}\to
C_{\textup{prod}}'=C_{\textup{prod}}$. 
Thus there is a resulting phase $\nu=C_{\textup{prod}}/|C_{\textup{prod}}|$.

\begin{lemma}
If $\Psi_1$ and $\Psi_2$ are $T^2$ representation with the same
$(R,n,k,\beta,\nu)$ then they are equivalent.

If $\Psi$ is a $T^2$ representation we can find an
equivalent representation $\tilde\Psi$ with the same 
$(R,n,k,\nu)$ and with 
$\beta$ replaced by $\tilde\beta$ such that 
$\pi-2\pi/n<\tilde\beta-\tfrac12\alpha\le\pi$.
\end{lemma}

\begin{proof}
\label{lm_T2_basic}
Since $(R,n,k,\nu)$ are the same for $\Psi_1$ and
$\Psi_2$ then $\alpha=2\pi k/n$ is the same for $\Psi_1$ and
$\Psi_2$.  By the same argument in lemma \ref{lm_S2_basic} we can
choose basis elements (\ref{rep_kets}) such that
$\CC(\beta+m\alpha)\in\Real_+$ for $m=1\sdotsc n-1$ and thus
$\CC(\beta)=\nu|\CC(\beta)|$. Doing this for both $\Psi_1$ and
$\Psi_2$ then the action (\ref{rep_u}), (\ref{rep_a_pm}) on these
basis elements are the same, therefore the representations are
equivalent.

In the proof of theorem \ref{th_reps} in the case of a $T^2$
representation when all $\CC(\theta)\ne0$ we choose $\beta$ to be any
value such that $e^{i\beta}$ was an eigenvalue of $\Psi(\Au)$. Since
the $n$ roots are equally spaced around the circle we can choose any
arc of length $2\pi/n$ to place $\beta$ in.  The one chosen makes the
calculations below simpler.
\end{proof}

For this subsection we assume we are given $(R,n,k,\nu)$, and
we wish to find $\beta$ so that $\Psi(R,n,k,\beta,\nu)$
is $T^2$-type irreducible representation. Again we define
$\beta'=\beta-\tfrac12\alpha$.

\begin{lemma}
Given $R$, $n$ and $k$ such that $R<\cos(\pi/n)\sec(\pi k/n)$ there
are no $T^2$ representation $\Psi(R,n,k,\beta,\nu)$.

Given $R$, $n$ and $k$ such that 
$\cos(\pi/n)\sec(\pi k/n)<R\le\sec(\pi k/n)$
then there exist a one parameter set of $T^2$ representation
$\Psi(R,n,k,\beta,\nu)$ with $\beta$ in the range
%[
\begin{align}
\pi-2\pi/n+\tfrac12\delta<\beta'<\pi-\tfrac12\delta
\sp
\textup{ where }
\cos(\tfrac12\delta)=R\cos(\pi k/n)
\label{rep_T2_beta}
\end{align}
%]

Given $R$, $n$ and $k$ such that $R>\sec(\pi k/n)$ there exist a
parameter set of irreducible $T^2$ representation
$\Psi(R,n,k,\beta,\nu)$ with the full range of $\beta$,
i.e. $\pi-2\pi/n<\beta'\le\pi$.
\end{lemma}
\begin{proof}
By looking at (\ref{rep_C}) we see that $\Psi(R,n,k,\beta,\nu)$ is an
irreducible $T^2$ representation if only if
%[
\begin{align}
\cos(\beta'+2\pi mk/n) + R\cos(\pi k/n) > 0
\sp
m=0\sdotsc n-1
\label{rep_T2_ineq}
\end{align}
%]
hence if $R>\sec(\pi k/n)$, (\ref{rep_T2_ineq}) is satisfied for 
all $\beta'$, hence $\Psi(R,n,k,\beta,\nu)$ is a representation for
all $\beta$.

For $R\le\sec(\pi k/n)$ then
similar to lemma \ref{lm_S2_non_min_R} if we set
$\theta'=(\beta'+2\pi mk/n)\!\!\!\!\mod 2\pi$ then $\theta'$ must lie either
in the range 
$0\le\theta'<\pi-\tfrac12\delta$ or $\pi+\tfrac12\delta<\theta'<2\pi$
where $\delta$ is given in
(\ref{rep_T2_beta}).  Since the set $\Set{e^{i(\beta'+2\pi
mk/n)}|m=0\sdotsc n=1}$ are equally spaced then $\delta<
2\pi/n$. Hence $0<\tfrac12\delta<\pi/n\le\tfrac12\alpha<\tfrac12\pi$
and $\cos$ is strictly decreasing we have
$0<\cos(\pi/n)<\cos(\tfrac12\delta)<1$. Since
$\cos(\tfrac12\delta)=R\cos(\pi k/n)$ we have $R>\cos(\pi/n)\sec(\pi
k/n)$.

Since $\beta'\le\pi$ then $\beta'\le\pi-\tfrac12\delta$. 
Now there exists and $m$ such that
$e^{i(\beta'+2\pi m k/n)}=e^{i(\beta'+2\pi/n)}$, hence
$(\beta'+2\pi mk/n)\!\!\!\! \mod
2\pi=\beta'+2\pi/n$ hence $\beta'+2\pi/n$ must lie either in the range
$0\le\beta'+2\pi/n<\pi-\tfrac12\delta$ or the range 
$\pi+\tfrac12\delta<\beta'+2\pi/n<2\pi$ but since $\beta'+2\pi/n>\pi$
we have $\beta'>\pi+\tfrac12\delta+2\pi/n$. Thus for $\Psi$ to be $T^2$
representation we must have $\beta'$ given by (\ref{rep_T2_beta}).

If $\beta'$ is given by (\ref{rep_T2_beta}) then all the $\theta'$ lie
in the permitted regions hence (\ref{rep_T2_ineq}) so it is a
representation. 
\end{proof}

We can see that $\Psi$ is a minimal $T^2$ representation if
$k=1$. Thus we have the following corollary
\begin{corrol}
If  $R\le1$ there are no $T^2$ representation. 

Given $R$ and $n$ such that $1<R\le\sec(\pi/n)$ then there exist a one
parameter set of irreducible $T^2$ minimal representation
$\Psi(R,n,k=1,\beta,\nu)$ with $\beta$ given in (\ref{rep_T2_beta}).

Given $R$ and $n$ such that $R>\sec(\pi/n)$ there exist a one
parameter set of irreducible $T^2$ representation
$\Psi(R,n,k=1,\beta,\nu)$ with the full range of $\beta$,
i.e. $\pi-2\pi/n<\beta'\le\pi$.
\end{corrol}
The range of possible $(R,\alpha)$ for minimal $T^2$ representation is
simply $R>1$ and $\alpha=2\pi/n$. These are pictures in figure
\ref{fig_T2}a. Figure \ref{fig_T2}b gives the range of possible
$(R,\alpha)$ for $n=11$. Figure \ref{fig_T2}c gives an example of a
$T^2$ representation.

%%%%%%%%%%%%%%%%%%%%%%%%%%%%%%%%%%%%%%%%%%%%%%%%%%

\subsection{Infinite dimensional $T^2$ representations}

Like the noncommutative torus $\Alg(R)$ also have infinite dimensional
representations.
We shall consider here only the infinite dimensional
representations of the form (with $\Hil$ having the basis
$\Set{\ket(\beta+m\alpha)\,|\,m\in\Intg}$)
%[
\begin{equation}
\begin{aligned}
&\Au\ket(\beta+m\alpha)=e^{i(\beta+m\alpha)}\ket(\beta+m\alpha) 
\sp
\Aa_-\ket(\beta+m\alpha) =
\CC(\beta+m\alpha)\ket({\beta+(m-1)\alpha})
\sp
\\&
\Aa_+\ket(\beta+m\alpha) =
\CC({\beta+(m+1)\alpha})\ket({\beta+(m+1)\alpha})
\end{aligned}
\label{rep_inf}
\end{equation}
%]
where 
%[
\begin{align}
C_{\beta+m\alpha} &= 
(\sec(\tfrac12\alpha)\cos(\beta'+m\alpha) + R)^{1/2}
\sp
\beta'=\beta-\tfrac12\alpha
\end{align}
%]
and where $\alpha\ne2\pi n/k$ for any $n,k\in\Intg$.
The eigenvalues of $\Psi(\Au)$ are dense on the circle, so there exist
an infinite dimensional representations if only if
$R\ge\sec(\alpha)$.

If $R=\sec(\alpha)$ there also exist semi-infinite
dimensional $T^2$ representation. These are given by $\beta'=-\pi$. Thus
$\Aa_-\ket(\beta)=0$ and $\Aa_+\ket(\beta-\alpha)=0$ so we can reduce the
Hilbert space to the subspaces
$\textup{span}\Set{\ket(-\pi+r\alpha),r\ge0}$ and 
$\textup{span}\Set{\ket(-\pi+r\alpha),r\le-1}$. 

%%%%%%%%%%%%%%%%%%%%%%%%%%%%%%%%%%%%%%%%%%%%%%%%%%%%%%%%%%%%%%%%%%%%%%

\newcounter{fff}
\def\footn{{\stepcounter{fff}${}^{\alph{fff}}$}}
\def\footnn{{${}^{\alph{fff}}$}}

\newlength{\tlen}
\newlength{\clen}
\setlength{\tlen}{3.5cm}
\setlength{\clen}{1.6cm}
\setlength{\unitlength}{1em}

\def\Qno{{\LARGE $\chi$}}

\def\Qyes{\makebox[3em][l]{\LARGE \raisebox{-0.2em}{
           {\raisebox{-0.4em}{$\ddot\smile\!\!\!$}\circle{2}}}}}
\def\Re{{R_\varepsilon}}

\begin{table}[t]
\begin{minipage}{\textwidth}
\begin{tabular}{|p{\tlen}|p{\clen}|p{\clen}|p{\clen}|p{\clen}|p{\clen}|p{\clen}|p{\clen}|}
\hline
&
Null &
Point &
Sphere &
Variety &
\multi{Sphere-\\Torus} &
\multi{Sphere-\\Torus\\bdd} &
Torus 
\\
\hline
\multi{$R$ range\\\ $\Re=(1\!+\!\varepsilon^2)^{1/2}$} &
{\footnotesize{$R\le-1$}} &
{\footnotesize{$R=-1$}} &
{\footnotesize{$-1\!<\!R\!<\!1$}} & 
{\footnotesize{$R=1$}} &
\multi{\footnotesize$1<R<\Re$} &
\multi{\footnotesize$R\!=\Re$} &
\multi{\footnotesize$R\!>\Re$}
\\
\hline
\multi{Topology of\\manifold $\man(R)$}
&
Null\footn
&
Point
&
$S^2$ top\footn
&
\multi{2dim\\variety\footn}
&
-\footn
&
-\footnn
&
$T^2$
\\
\hline
\multi{Minimal $S^2$\\representations.}
& \Qno
& \Qno
& \Qyes
& \Qyes
& \Qyes
& \Qno
& \Qno
\\
\hline
\multi{Non minimal $S^2$\\representations.}
& \Qno
& \Qno
& \Qno
& \Qno
& \Qyes
& \Qyes
& \Qno
\\
\hline
\multi{Finite dim  $T^2$\\ representations}
& \Qno
& \Qno
& \Qno
& \Qno
& \Qyes\footn
& \Qyes\footn
& \Qyes
\\
\hline
\multi{Semi $\infty$-dim $T^2$ \\representations}
& \Qno
& \Qno
& \Qno
& \Qno
& \Qno
& \Qyes
& \Qno
\\
\hline
\multi{$\infty$-dim $T^2$ \\representations}
& \Qno
& \Qno
& \Qno
& \Qno
& \Qno
& \Qyes
& \Qyes
\\
\hline
\end{tabular}
\footnotesize{ 
\begin{tabular}{lll}
Notes: & {a: No solution to equation} 
& {b: $\man(R)$ is convex if $R<0$.} \\
& {c: The point $(x,y,z)=(0,0,0)$ not smooth.} 
& {d: Since $(1+\varepsilon^2)^{1/2}=1$ so no sphere-torus exists.} \\
& {e: $\beta$ range is limited.} 
& {f: $\beta$ range excludes one point.}
\end{tabular}
}
\end{minipage}
\caption{Summary of possible representations for different $R$ and
$\varepsilon$.} 
\label{fig_concl}
\end{table}

\section{Conclusion and Discussion}
\label{ch_concl}
Table \ref{fig_concl} gives a list of all the possible
representations. 
We can see that representations reflect the topology, but not
completely. 

There are loosely speaking four regions of $R$.
If $R\le-1$ then there is either no manifold $\man(R)$ or
it is just a point. Consequently there are no representations either.

The next region $-1<R<1$ the algebra $\Alg(R)$ is closest to the
noncommutative sphere. The commutative limit is the sphere, and
there exist only minimal $S^2$ representations
which are the closest to the representations of the noncommutative
sphere, in that they are parameterised by $n$. 

In the region $1<R<(1+\varepsilon^2)^{1/2}$ the algebra $\Alg(R)$ is a
new object which we may call the ``sphere-torus''. This is a purely
noncommutative object since setting $\varepsilon=0$ gives $1<R<1$ so
there is no commutative analogue. In this region $\Alg(R)$ has
the minimal $S^2$ representations, the non minimal $S^2$
representations and the finite dimensional $T^2$ representations. 

Only when $R>(1+\varepsilon^2)^{1/2}$ can we say that
$\Alg(R)$ is a deformation of the torus. There are no $S^2$
representations and, as well as the finite dimensional $T^2$
representations, the infinite $T^2$ representations exist.

If any definition of genus is applied to noncommutative geometry it
would be interesting to see what values it would attain for the
sphere-torus.

%%%%%%%%%%%%%%%%%%%%%%%%%%%%%%%%%%%%%%%%%%%%%%%%%%%%%%%%%%%%%%%%%%%%%%
\section*{Acknowledgements}

The author would like to thank Arne Sletsj{\o}e Math Department 
University of Oslo, John Madore Laboratoire de Physique Th{\'e}orique 
Universit{\'e} de Paris Sud and Richard Kerner Laboratoire
Physique Th{\'e}orique des Liquides, Universit{\'e} Pierre et Marie
Curie, Paris VI for help in preparing this paper. 

The author would like to thank the the math department of the 
University of Oslo for their facilities.

%%%%%%%%%%%%%%%%%%%%%%%%%%%%%%%%%%%%%%%%%%%%%%%%%%%%%%%%%%%%%%%%%%%%%%


\begin{thebibliography}{9}

\bibitem{Madore_book} J. Madore 1995, 
{\it An Introduction to Noncommutative Differential Geometry and its
Physical Applications},
 Cambridge University Press.

\bibitem{Sternheimer1} D. Sternheimer,
{\it Deformation Quatization: Twenty Years After}
math.QA/9809056 (And references therein.)

\bibitem{Gratus1} J. Gratus, 
{\it A Natural Basis of States for the Noncommutative Sphere 
and its Moyal bracket} 
J. Maths. Phys {\bf 38} (8), (Aug 1997) 4283 - 4300, 
q-alg/9703038

\end{thebibliography}
\end{document}